\documentclass{amsart}

\usepackage{amssymb}

\newtheorem{thm}{Theorem}[section]
\newtheorem{lem}[thm]{ Lemma}
\newtheorem{lemma}[thm]{Lemma}
\newtheorem{cor}[thm]{Corollary}
\newtheorem{prop}[thm]{Proposition}
\newtheorem{defi}[thm]{Definition}

\def\Re{{\rm Re}\,}
\def\Im{{\rm Im}\,}
\def\R{{\mathbb R}}
\def\calS{{\mathcal S}}
\def\cH{{\mathcal H}}
\def\pr{\partial}
\def\la{\langle}
\def\ra{\rangle}
\def\les{\lesssim}
\def\C{{\mathbb C}}
\def\eps{\varepsilon}
\def\gtr{\gtrsim}
\def\supp{{\rm supp}}
\def\const{{\rm const}}
\def\wt{\widetilde}
\def\calM{{\mathcal M}}
\def\calH{{\mathcal H}}
\def\walze{-\frac{d}{2}-\sigma}
\def\walz{-\frac{d+1}{2}-\sigma}
\def\calL{{\mathcal L}}

\begin{document}
\numberwithin{equation}{section}

\title[Decay for the Wave and Schr\"odinger evolutions: part II]{Decay
for the wave and Schr\"odinger evolutions on manifolds with conical ends, Part~II}

\author{Wilhelm Schlag}
\address{University of Chicago, Department of Mathematics,
5734 South University Avenue, Chicago, IL 60637, U.S.A.}
\email{schlag@math.uchicago.edu}
\thanks{The first author was partly supported by the National
Science Foundation DMS-0617854.}

\author{Avy Soffer}
\address{Rutgers University, Department of Mathematics, 110 Freylinghuysen Road, Piscataway, NJ 08854, U.S.A.}
\email{soffer@math.rutgers.edu}
\thanks{The second author was partly supported by the National
Science Foundation DMS-0501043.}

\author{Wolfgang Staubach}
\address{Department of Mathematics,
Colin Maclaurin Building, Heriot-Watt University, Edinburgh, EH14 4AS}
\email{W.Staubach@hw.ac.uk}


\date{}

\begin{abstract}
Let $\Omega\subset \R^N$ be a compact imbedded Riemannian manifold of
dimension~$d\ge1$ and define the $(d+1)$-dimensional Riemannian manifold
$\calM:=\{(x,r(x)\omega)\::\: x\in\R,\, \omega\in\Omega\}$ with  $r>0$ and smooth, and the natural metric
$ds^2=(1+r'(x)^2)dx^2+r^2(x)ds_\Omega^2$.  We require that $\calM$ has
conical ends: $r(x)=|x| + O(x^{-1})$ as $x\to \pm\infty$. The Hamiltonian flow on such manifolds
always exhibits trapping.  Dispersive
estimates for the Schr\"odinger evolution $e^{it\Delta_\calM}$ and
the wave evolution $e^{it\sqrt{-\Delta_\calM}}$ are obtained for
data of the form $f(x,\omega)=Y_n(\omega) u(x)$ where $Y_n$ are eigenfunctions of~$-\Delta_\Omega$ with eigenvalues~$\mu_n^2$.
In this paper we discuss all
cases $d+n>1$. If $n\ne0$ there is the following {\em accelerated local decay} estimate: with
\[
0< \sigma = \sqrt{2\mu_n^2+(d-1)^2/4}-\frac{d-1}{2}
\]
and all $t\ge1$,
\[\Vert w_\sigma\,
e^{it\Delta_{\calM}}\, Y_nf\Vert_{L^{\infty}(\calM )} \le C(n,\calM,\sigma)\,
t^{-\frac{d+1}{2}-\sigma}\Vert w_\sigma^{-1} f\Vert_{L^1(\calM )}
\]
where $w_\sigma(x)=\la x\ra^{-\sigma}$, and similarly for the wave evolution.
Our method combines two main ingredients:

\noindent (A) a detailed scattering analysis of Schr\"odinger operators
of the form $-\partial_\xi^2 + (\nu^2-\frac14)\la\xi\ra^{-2}+U(\xi)$ on the line where $U$ is real-valued and smooth with $U^{(\ell)}(\xi)=O(\xi^{-3-\ell})$ for all
$\ell\ge0$
as $\xi\to\pm\infty$ and $\nu>0$. In particular, we introduce the notion of a zero energy resonance
for this class and derive an asymptotic expansion of the Wronskian between the outgoing Jost solutions as the energy tends to zero.
In particular, the division into Part~I and Part~II can be explained by the former being resonant at zero energy, where
the present paper deals with the nonresonant case.

\noindent (B) estimation of oscillatory integrals  by (non)stationary phase.
\end{abstract}

\maketitle

\section{Introduction}

As in Part~I, see~\cite{SSS}, we consider the following class of manifolds~$\calM$:

\begin{defi}
\label{def:cones} Let $\Omega\subset \R^N$
be an imbedded compact $d$-dimensional Riemannian manifold with metric $ds_\Omega^2$ and define the $(d+1)$-dimensional manifold
\[
\calM:=\{(x,r(x)\omega)\:|\: x\in\R,\; \omega\in\Omega\},\quad
ds^2=r^2(x)ds_{\Omega}^2 +(1+r'(x)^2)dx^2
\]
where $r\in C^\infty(\R)$ and $\inf_x r(x)>0$. We say that there is
a {\em conical end} at the right (or left) if
 \begin{equation}\label{eq:cone} r(x)=|x|\,(1+ h(x)),\quad
  h^{(k)}(x) =O(x^{-2-k}) \quad \forall\; k\geq 0
  \end{equation}
as $x\to\infty \; (x\to-\infty)$.
\end{defi}

Of course we can consider cones with arbitrary opening angles here but this adds
nothing of substance. Examples of such manifolds are given by surfaces of revolution with $\Omega=S^1$
such as the one-sheeted hyperboloid. They have the property that
the entire Hamiltonian flow on~$\calM$ is trapped on the
set $(x_0,r(x_0)\Omega)$ when $r'(x_0)=0$.

The main results of this paper are global in time dispersive
estimates for the Schr\"odinger evolution $e^{it\Delta_\calM}$ and
the wave evolution $e^{it\sqrt{-\Delta_\calM}}$, where
$\Delta_{\calM}$ denotes the Laplace-Beltrami operator on~$\calM$.
These results should be contrasted to the large number of papers
studying wave evolution on curved back grounds, see for example,
\cite{Bou}, \cite{BGT1}, \cite{BGT2}, \cite{CKS}, \cite{Doi},
\cite{G}, \cite{HTW1}, \cite{HTW2}, \cite{RZ},\cite{RodTao},
\cite{SS}, \cite{ST}, \cite{T}, \cite{DT}. However, these references
either consider the evolution on general manifolds for short times,
or the global evolution on asymptotically flat perturbations of the
Euclidean metric under a non-trapping condition.
However, see the recent papers \cite{NonZwo}, \cite{Chris1}--\cite{Chris3}, as well
as the more classical paper \cite{GerSjo} for semiclassical results close to a hyperbolic orbit
of the Hamiltonian flow, and also \cite{Burq} for local decay of the energy for the wave equation
without any assumption on trapping.

The main point here
is to carefully examine the long time behavior for a  class of
examples that do exhibit trapping. In Part~I we proved the case $d=1, n=0$ which is special
(it can be viewed as an end-point case of the theory in Section~\ref{sec:scatternu} below).

In what follows, $\{Y_n,\mu_n\}_{n=0}^\infty $
denote the $L^2$-normalized eigenfunctions and eigenvalues, respectively,  of~$\Delta_\Omega$.
In other words, $-\Delta_\Omega Y_n = \mu_n^2 Y_n$ where $0=\mu_0^2< \mu_1^2\le \mu_2^2\le\ldots$

\begin{thm}\label{thm1}
Let $\calM$ be asymptotically conical at both ends  in the sense of
Definition~\ref{def:cones} with $d\ge1$ arbitrary.  For each $d\ge1$, $n\ge0$, let
\[
\nu=\nu(d,n):= \sqrt{2\mu_n^2 +(d-1)^2/4}.
\]
For each  $n\ge0$ and all $0\le\sigma\le \nu(d,n) -\frac{d-1}{2}$,
there exist constants $C(n,\calM,\sigma)$ and $C_1(n,\calM,\sigma)$ such that for all $t>0$
\begin{align} \Vert w_\sigma\,
e^{it\Delta_{\calM}}\, Y_nf\Vert_{L^{\infty}(\calM )} &\le \frac{C(n,\calM,\sigma)}{t^{\frac{d+1}{2}+\sigma}}\Big\Vert \frac{f}{w_\sigma }\Big\Vert_{L^1(\calM )} \label{eq:schr_d}\\
\Vert w_\sigma\, e^{\pm it\sqrt{-\Delta_{\calM}}}\, Y_nf\Vert_{L^{\infty}(\calM )} &\le \frac{C_1(n,\calM,\sigma)}{t^{\frac{d}{2}+\sigma}}\Big(\Big\Vert \frac{f'}{w_\sigma }\Big\Vert_{L^1(\calM )} + \Big\Vert \frac{f}{w_\sigma }
\Big\Vert_{L^1(\calM)}  \Big) \label{eq:wave_d}
\end{align}
provided $f=f(x)$ does not depend on~$\omega$. Here $w_\sigma(x):=\la x\ra^{-\sigma}$ are weights on~$\calM$.
\end{thm}

In our previous paper~\cite{SSS} we dealt with the case $d=1, n=0$ and
proved~\eqref{eq:schr_d} and~\eqref{eq:wave_d} for that case. Needless to say, it is the analogue
of the usual dispersive decay estimate for the Schr\"odinger and wave
evolutions on~$\R^2$. Clearly, the {\bf local} decay given by $\sigma>0$ has
no analogue in the Euclidean setting and it also has no meaning for
$n=0$. To motivate it, one can try to rely on the geodesic flow
on~$\calM$.
 As an example, take~$\calM$ to be the one-sheeted hyperboloid. It has a unique
closed geodesic $\gamma_0$ at the neck and any other
geodesic~$\gamma$ that crosses $\gamma_0$ will pull away from it and
never return. The analogue of~$\mu_n=n$  would be the
velocity of~$\gamma$ and thus, the larger $n$, the faster $\gamma$
will pull away. This can serve to ``explain'' the improvement in terms of the power of~$t$
for $\sigma>0$ in the
following sense: imagine $f=f(x)$ to be a highly localized bump
function centered around $\gamma_0$. Then the Schr\"odinger flow
will (in phase space) mimic the geodesic flow on the cotangent
bundle, at least for short times. By the ``pulling away'' logic we
expect such data to disintegrate under the Schr\"odinger flow ---
and more strongly so as the angular momentum~$n$ increases. Hence,
as long as we only ask for the size of the solution close
to~$\gamma_0$ (this is the effect of the weight~$w_\sigma$) we would
expect to see very little of the wave left around the neck. Clearly,
this is a very much a negative curvature effect which should be contrasted
to $\calM=S^{d+1}$, for example.

However, this heuristic reasoning has to be taken with a grain of
salt. First --- and to the best of the authors' knowledge --- it is
 not clear how  to derive the exact power law above
  via a classical approximation. This is due to the
 fact that dispersive effects limit the accuracy of any classical
 approximation for long times (in more technical terms, this is the
 problem of constructing global parametrices by semi-classical
 methods).
 Second, note that Theorem~\ref{thm1} does not specify the behavior
 of~$\calM$ close to $x=0$. For example,
 Theorem~\ref{thm1} applies to a surface which is obtained as
 follows: cut the one-sheeted hyperboloid at the neck and glue the
 two pieces smoothly to a large sphere from which we have removed
 caps around the poles. The sphere of course has a continuum of
 stable closed geodesics. Nevertheless, due to the dispersion of the
 Schr\"odinger and wave flows the solution will spread into regions of~$\calS$
 that exhibit the aforementioned instability of the geodesic flow
 typical of negatively curved surfaces. Theorem~\ref{thm1} states
 that over {\em long times} the {\em power law is universal} and does not
 see the local geometry. On a more technical level, let us mention
 that~\eqref{eq:schr_d} and~\eqref{eq:wave_d} are optimal  with regard to both their respective $t^{-\frac{d+1}{2}-\sigma}$
and  $t^{-\frac{d}{2}-\sigma}$ decay rates,
 as well as the polynomial weights $w_\sigma$ (in the sense that we cannot
 choose a smaller power) and the range of allowed~$\sigma$.

 \noindent A subtle point arises here, which is the size of the constants
 $C(n,\calS,\sigma)$ and~$C_1(n,\calS,\sigma)$, especially with regard to their
 asymptotic behavior as $n\to\infty$.
We do not address this asymptotic issue in~$n$
 at all in this paper. In fact, the methods of this paper
 were not designed with a view towards optimal constants --- but rather to exhibit the correct asymptotic
 behavior in~$t$ --- and the constants $C(n,\calS,\sigma)$,  $C_1(n,\calS,\sigma)$ produced by our proof grow
 super-exponentially in~$n$. A forthcoming paper will address the question of how our constants depend on~$n$
by considering
 $\hbar=n^{-1}$ as a {\em small} semi-classical
 parameter.
Hence, it is appropriate to view this paper as dealing
 with the {\em intermediate} regime of $n$, namely those that are
 not zero but not too large.

We believe that
 the  analysis  of the Laplacean $\Delta_\calM$ which is carried out in Parts~I and~II
should be of independent interest and  the Schr\"odinger and wave flow merely serve as an
 example where our asymptotic analysis applies. Note in particular that in Section~\ref{sec:scatternu} we
develop the scattering theory of the class of Schr\"odinger operators on the line
\[
\calH_\nu = -\partial_\xi^2 + (\nu^2-\frac14)\la\xi\ra^{-2}-U_\nu(\xi), \qquad \frac{d^\ell U_\nu(\xi)}{d\xi^\ell} =O(\xi^{-3-\ell})\] for all
$\ell\ge0$
as $\xi\to\pm\infty$ and $\nu>0$ (the $\nu$ is defined in Theorem~\ref{thm1}). Section~\ref{sec:scatternu} is ``abstract'' in the sense
that it does not draw on anything from other sections. We obtain approximations to the Jost solutions of $\calH_\nu$
as the energy $\lambda$ tends to zero and also find that their Wronskian is of the form $\lambda^{1-2\nu}$ {\em provided
 there is no zero energy resonance}. By this we mean that the two subordinate solutions of the equation $\calH_\nu f=0$
as $\xi\to\pm\infty$, respectively, do not form a {\em globally} subordinate solution on~$\R$. See Definition~\ref{def:resonance} below.
Although the notion of a zero energy resonance is standard for potentials that belong to $\{\la\xi\ra^{-1}V(\xi)\:|\:V\in L^1\}$, we are
not aware of a reference where the conclusions of Section~\ref{sec:scatternu} are reached. This notion also helps to
explain the difference between Part~I and Part~II: the former is {\em resonant} whereas the latter deals with the {\em nonresonant case}.
This is in agreement with the fact that the Laplacian on $\R^2$ has a zero energy resonance, whereas on $\R^n$ with $n\ge3$ it does not.
The fact that Part~I is {\em resonant}  is due to the fact that the zero energy solutions are $u_0(\xi)=\sqrt{r(\xi)}$,
and $u_1(\xi)=\sqrt{r(\xi)}\int_0^\xi \frac{d\eta}{r(\eta)}$. The subordinate one is $u_0$ which is global.

In the context of our conical manifolds we are able to settle the important resonant vs.\ nonresonant  question by knowledge of the {\em zero energy solutions}
of $\calH_\nu$ which of course is equivalent to knowledge of the (spherical) harmonics of $\Delta_\calM$. In fact, we shall see later on that for all $\nu>0$
the maximum principle allows us to conclude that we are in nonresonant case (indeed, in the resonant case there would need to be a nonzero harmonic
function on~$\calM$ that vanishes at both ends which contradicts the maximum principle). As an example, for $d=1, n>0$
these functions are (with $\mu_n=n$ since $\Omega$ is isometric to $S^1$ for $d=1$)
\[
\calH_{1,n} (r^{\frac12} e^{\pm n y})=0,\quad y(\xi)=\int_0^\xi \frac{d\eta}{r(\eta)}
\]
Because $y$ is odd, the smaller branch at $\xi=\infty$
has to be larger one at $\xi=-\infty$ which places us in the nonresonant case.

Finally, let us remark that the methods of this paper cannot touch
 non-rotationally invariant perturbations of the metric on surfaces of
 revolution, let alone a non-symmetric example like three
 half-cones glued together smoothly (a ``conical three-foil'').
Another example would be two parallel planes joined by $k$ necks.
While the case $k=1$ is essentially covered by Theorem~\ref{thm1},
the cases $k \ge2$ are of course very different. It would be most interesting to find a way of
 approaching an analogue of Theorem~\ref{thm1} for  manifolds which do not allow for separation
of variables as we use here.

\section{The setup and an overview over the method}
\label{sec:setup}

For the convenience of the reader, we reproduce some of the material from Section~2 of Part~I.
First, recall that the Laplace-Beltrami operator on $\calM$ where the base
$\Omega$ is of dimension $d\ge1$, is

\begin{equation}\label{eq:LapBelt}
\Delta_{\calM}=\frac{1}{r^d(x)\sqrt{1+r'(x)^2}}\,
\partial_x\left(\frac{r^d(x)}{\sqrt{1+r'(x)^2}}
\partial_x\right)+\frac{1}{r^2(x)}\Delta_\Omega .
\end{equation} In arclength parametrization
$$\xi (x)=\int_0^x\sqrt{1+r'(y)^2}\, dy$$
\noindent  (\ref{eq:LapBelt}) reads
\begin{equation}\label{eq:LapBelt2}
\Delta_{\calM}=\frac{1}{r^d(\xi)}\partial_{\xi}(r^d(\xi
)\partial_{\xi})+\frac{1}{r^2(\xi )}\Delta_\Omega \end{equation} \noindent
where we have abused notation: $r(\xi )$ instead of $r(x(\xi ))$.
Setting $\rho (\xi ):=\frac{d}{2}\frac{\dot{r}(\xi )}{r(\xi )}$
yields \begin{equation}\label{2.1}\Delta_{\calM}\,y(\xi ,\omega )=
\partial^2_{\xi}y+2\rho\partial_{\xi}y+\frac{1}{r^2}\Delta_\Omega y.\end{equation}
The first order term in \eqref{2.1} is removed by setting
\begin{equation}\label{3} y(\xi ,\omega )=r(\xi )^{-\frac{d}{2}}u(\xi ,\omega).
\end{equation}
 Then \begin{equation}\label{eq:4} \Delta_\calM y =
\partial^2_{\xi}y+2\rho\partial_{\xi}y+\frac{1}{r^2}\Delta_{\Omega} y=r^{-d/2}[-\calH u+\frac{1}{r^2}\Delta_\Omega u]
\end{equation} \noindent with
\begin{equation}\label{eq:calH} V_1(\xi ):=\rho^2(\xi
)+\dot{\rho}(\xi ),\quad \calH=-\partial^2_{\xi}+V_1. \end{equation} Note that
the Schr\"odinger operator $\calH$ can be factorized as
\begin{equation}\label{eq:calL}
\calH =\calL^*\calL,\quad \calL = -\frac{d}{d\xi} + \rho
\end{equation}
In particular, $\calH$ has no negative spectrum.  Now specialize further to $u(\xi,\omega)=Y_n(\omega) \phi(\xi)$.
Then
\[
\calH u-\frac{1}{r^2}\Delta_\Omega u = Y_n \calH_{d,n} \phi, \qquad \calH_{d,n} = -\partial_\xi^2 + V, \quad V(\xi):= V_1(\xi)+\frac{\mu_n^2}{r^2(\xi)}.
\]
The Schr\"odinger operator $\calH_{d,n}$ is of fundamental importance to this paper. It has a smooth potential $V$
with the following asymptotic behavior:
In Part~I we proved that, see Lemma~2.2 and Corollary~2.3 of~\cite{SSS},
\[
V(\xi)=V_1(\xi) + \frac{\mu_n^2}{r^2(\xi)} = (2\mu_n^2 + d(d-2)/4)\la \xi\ra^{-2} + O(\la\xi\ra^{-3}) = \big(\nu^2-\frac14\big)\la\xi\ra^{-2} + O(\la\xi\ra^{-3}).
\]
Here $\nu^2 = 2\mu_n^2 + (d-1)^2/4$ is exactly as in Theorem~\ref{thm1} and the $O(\cdot)$ term behaves like a symbol, which means that
\[
|\partial_\xi^\ell O(\la\xi\ra^{-3})| \le C\la\xi\ra^{-3-\ell},\quad \forall\;\ell\ge0
\]
Note carefully that $d+n>1$ corresponds precisely to $\nu>0$ in~$\calH_{d,n}$.
In terms of the
Schr\"odinger evolution,
\[
e^{-it\Delta_\calM}\, Y_n f = r^{-\frac{d}{2}} Y_n\,  e^{it \calH_{d,n}}\, r^{\frac{d}{2}}
f \quad \forall\; f=f(\xi)
\]
and similarly for the wave equation. In particular, any estimate of
the form
\[
\big\| w_\sigma\, e^{-it\Delta_\calM} Y_n f \|_{L^\infty(\calM)} \le Ct^{-\alpha}
\| \frac{f}{w_\sigma }\|_{L^1(\calM)} \quad \forall\; t>0,\; f=f(\xi)
\] with arbitrary $\alpha\ge0$ and some constant $C$ that does not depend on~$t$,
is equivalent to one of the form
\begin{equation}\label{11}
\big\|r^{-\frac{d}{2}} w_\sigma e^{it\calH_{d,n}}\,r^{-\frac{d}{2}} u\big\|_{L^{\infty}(\R
)} \le C'\, t^{-\alpha} \|\frac{u}{w_\sigma }\|_{L^1(\R)} \quad \forall\; t>0,\;
u=u(\xi)\end{equation} with a possibly different constant $C'$. Here we
absorbed the weight from the volume element $dv_\calM = r^d d\xi
dv_\Omega$ arising in the $L^1(\calM)$ norm into the left-hand side
of~\eqref{11}. An analogous reduction is of course valid for the
wave evolution.  As usual, the functional calculus applied
to~\eqref{11} yields
\[
e^{it\calH_{d,n}} = \int_0^\infty e^{it\lambda} E(d\lambda)
\]
where $E(d\lambda)$ is the spectral resolution of~$\calH_{d,n}$. The point
is that there is an ``explicit expression'' for~$E(d\lambda)$:
\[
E(d\lambda^2)(\xi,\xi') = 2\lambda \Big\{
\Im\Big[\frac{f_{+,\nu}(\xi,\lambda)f_{-,\nu}(\xi',\lambda )}{W_\nu(\lambda
)}\Big]\chi_{[\xi>\xi']} + \Im
\Big[\frac{f_{-,\nu}(\xi,\lambda)f_{+,\nu}(\xi',\lambda )}{W_\nu(\lambda )}\Big]
\chi_{[\xi<\xi']} \Big\}\,d\lambda
\]
where \[ W_\nu(\lambda ):=W(f_{-,\nu}(\cdot ,\lambda),f_{+,\nu}(\cdot ,\lambda))=
f_{+,\nu}'(\cdot ,\lambda) f_{-,\nu}(\cdot ,\lambda) - f_{-,\nu}'(\cdot ,\lambda)
f_{+,\nu}(\cdot ,\lambda)
\] is the Wronskian of the {\em Jost solutions} $f_{\pm,\nu}(\cdot ,\lambda )$ of
the following ordinary differential equation
\begin{equation}
\begin{aligned}\label{13}
\calH_{d,n}\, f_{\pm,\nu}(\xi ,\lambda )&=-f_{\pm,\nu}''(\xi ,\lambda )+V(\xi
)f_{\pm,\nu}(\xi ,\lambda ) =\lambda^2\, f_{\pm,\nu}(\xi ,\lambda)\\
f_{\pm,\nu}(\xi ,\lambda ) &\sim e^{\pm i\lambda\xi}\qquad \mathrm{as}\;
\xi\to\pm\infty
\end{aligned}
\end{equation}
 \noindent provided $\lambda\neq 0$. It is a standard fact that
these Jost solutions exist because of the decay
$|V(\xi )|\les\langle\xi\rangle^{-2}.$
In fact, they are easily seen to exist provided the perturbing potential $V$ is
in~$L^1$, see~\cite{DT}. Moreover, these Jost solutions are
continuous in the energy~$\lambda$ as $\lambda\to0$ precisely when
$\la\xi\ra V(\xi)\in L^1(\R)$ --- which obviously fails here. On the
other hand, it is common knowledge that the {\em asymptotic form},
i.e., the $t\to\infty$ decay law, of  dispersive estimates
like those in Theorem~\ref{thm1} crucially depend on the  behavior of
the spectral measure as $\lambda\to0$. These two facts of course fit
together well, as anything unusual like the accelerated local decay given by~$\sigma>0$ must be
reflected in the Jost solutions $f_{\pm,\nu}(\cdot,\lambda)$ around $\lambda=0$. In fact,
we show below that their Wronskian  for $\lambda>0$ and with
$\nu$ as in Theorem~\ref{thm1} satifies
\begin{equation}\label{eq:wronsk}
W(f_{+,\nu}(\cdot,\lambda), f_{-,\nu}(\cdot,\lambda)) = c_\nu\,
\lambda^{1-2\nu} (1+O(\lambda^\eps))\text{\ \ as\ \ }\lambda\to0+
\end{equation}
Here $\eps>0$ is small depending on~$\nu$,
cf.~Proposition~\ref{prop:wronski} and Corollary~\ref{cor:Wnu_nonres} below. An intuitive way of
viewing the leading order behavior in~\eqref{eq:wronsk} is as
follows. From elementary quantum mechanics considerations we expect this leading order  to be
given by $\lambda e^{d_A(\lambda)}$ where $d_A(\lambda)$ is the {\em
Agmon distance} between the {\em turning points} of~$\cH_{d,n}$ at
energy $\lambda^2$ (actually, it turns out that one needs to ignore the $-\frac14\xi^{-2}$ piece of $V$ for that purpose).
Recall that the turning points  $ \xi_1<0< \xi_2$
are determined from the relation  $\nu^2\la \xi_j\ra^{-2}=\lambda^2$ for $j=1,2$. The
Agmon distance between $\xi_1$ and $\xi_2$ is then defined to be
\[
d_A(\lambda) := \int_{\xi_1}^{\xi_2} \sqrt{\nu^2\la\xi\ra^{-2}
-\lambda^2}\, d\xi
\]
and thus
\[
d_A(\lambda) = 2\nu |\log \lambda|
\]
to leading order as $\lambda\to0$. Finally, this exactly gives
$\lambda e^{-2\nu\log\lambda}=\lambda^{1-2\nu}$ for the Wronskian as
claimed, see~\eqref{eq:wronsk}.
We caution the reader, though, that this heuristic via the Agmon distance only applies
in the nonresonant case. See Section~\ref{sec:scatternu}, in particular Definition~\ref{def:resonance}.

In view of the preceding, the estimates \eqref{eq:schr_d} and~\eqref{eq:wave_d} of Theorem~\ref{thm1} reduce to the following
respective oscillatory integral estimates\footnote{We remark that the imaginary part in \eqref{eq:crux} is crucial.}, uniformly in $\xi>\xi'$ (the case
$\xi\le \xi'$ be analogous)
\begin{align}
  \label{eq:crux} &\left|\int_0^\infty e^{it\lambda^2} \lambda\,  \Im \Big[
\frac{f_{+,\nu}(\xi,\lambda)
    f_{-,\nu}(\xi',\lambda)}{W_\nu(\lambda)} \Big]\,d\lambda   \right|
    \le C(\nu,\calM,\sigma)\, (\la\xi\ra\la\xi'\ra)^{\frac{d}{2}+\sigma} t^{\walz}\\
\label{eq:crux_wave} & \left|\int_{-\infty}^\xi \int_0^\infty e^{\pm it\lambda} \lambda\,  \Im \Big[
\frac{f_{+,\nu}(\xi,\lambda)
    f_{-,\nu}(\xi',\lambda)}{W_\nu(\lambda)} \Big]\,d\lambda\, (\la\xi\ra\la\xi'\ra)^{\walze} \, \phi(\xi')\, d\xi' \right|
    \\&\le C(\nu,\calM,\sigma)\,t^{\walze} \int (|\phi'(\eta)|+|\phi(\eta)|)\, d\eta\nonumber
\end{align}
It turns out that the  $t^{-\nu}$ improvement over the usual $t^{-\frac{d+1}{2}}$ decay
in~\eqref{eq:crux} stems from the $\lambda^{-2\nu}$ appearing in
the Wronskian~\eqref{eq:wronsk}; indeed, we prove in this paper that
the $\xi=\xi'=0$ case of~\eqref{eq:crux}  reduces to the standard stationary phase type bound
\[
\Big|\int_0^\infty e^{it\lambda^2} \lambda^{1+2\nu}
\chi(\lambda)\,d\lambda \Big|\le C t^{-1-\nu}
\]
where $\chi$ is a smooth cut-off function to the interval $[0,1]$,
say. Since $\frac{d+1}{2}+\sigma\le 1+\nu$, this estimate implies the desired $t^{\walz}$ bound from~\eqref{eq:schr_d}.
Note that this calculation also show the optimality of the upper bound $\sigma\le \nu-\frac{d-1}{2}$.

The reader should compare $W_\nu$ in~\eqref{eq:wronsk} with the
Wronskian for $n=0, d=1$ derived in~\cite{SSS}:
\[
W(\lambda )=2\lambda \left( 1+ic_3+i\frac{2}{\pi}\log\lambda
\right)+O(\lambda^{\frac{3}{2}-\eps}) \text{\ \ as\ \ }\lambda\to0+
\]
On a technical level, the logarithmic term in $\lambda$ makes the
$n=0, d=1$ case of~\eqref{eq:crux} somewhat harder to analyze than the cases $d+n>1$
(as already mentioned, $d+n=1$ is exactly $\nu=0$). Not surprisingly, in  proving dispersive estimates
for $-\Delta_{\R^2}+V$ one encounters similar logarithmic issues,
see~\cite{Sch}.

\section{The scattering theory of $\cH_\nu,\; \nu>0$}
\label{sec:scatternu}

This section can and should be
viewed as a separate entity, as it is kept completely general without any reference to the other sections.
Our goal is to develop the scattering theory of the following class of operators ($\nu=0$ is treated in Part~I, see~\cite{SSS}):

\begin{defi}\label{def:Hnu}  We define the class of operators
\[
 \calH_\nu := -\partial_\xi^2 + V(\xi)
\] where
 $V\in C^\infty(\R)$ is real-valued  with the property that
\[
V(\xi)= \big(\nu^2-\frac14\big)\xi^{-2} - U_\nu(\xi), \quad U_\nu\in C^\infty(\R\setminus\{0\})
\]
with $\nu>0$, $U_\nu$ real-valued and $U_\nu^{(\ell)}(\xi)=O(\xi^{-3-\ell})$ for all $\ell\ge0$  as $\xi\to\pm\infty$.
\end{defi}

The goal here is to obtain representations of the Jost solutions $f_{\pm,\nu}(\xi,\lambda)$ of $\calH_\nu$ and their Wronskian
$W_\nu(\lambda)$ especially as $\lambda\to0$.
Our main results concerning the class~$\calH_{\nu}$ are
Proposition~\ref{prop:wronski} and Corollary~\ref{cor:Wnu_nonres}.
We begin with certain bases of zero energy solutions $u_{0,\nu}^{\pm}$ and $u_{1,\nu}^{\pm}$.
All functions are smooth in $\xi$ where they are defined.

\begin{lemma}\label{lem:zero_en}
There are solutions $u_{0,\nu}^{\pm}$ and $u_{1,\nu}^{\pm}$ of $\calH_\nu f=0$ with the following properties:
\begin{align}
u_{0,\nu}^{+}(\xi) & = \xi^{\frac12+\nu}(1+O(\xi^{-\alpha})), \quad  u_{1,\nu}^{+}(\xi) = \xi^{\frac12-\nu}(1+O(\xi^{-1}))  \text{\ as\ }\xi\to\infty \\
u_{0,\nu}^{-}(\xi) & = |\xi|^{\frac12+\nu}(1+O(\xi^{-\alpha})), \quad  u_{1,\nu}^{-}(\xi) = |\xi|^{\frac12-\nu}(1+O(\xi^{-1}))  \text{\ as\ }\xi\to-\infty
\end{align}
The $O(\cdot)$ terms behave like symbols under differentiation in $\xi$ and $0<\alpha\le \min(2\nu,1)$. Furthermore, the solutions $u_{1,\nu}^\pm$ are unique
with the stated asymptotic behavior and
\[
W(u_{0,\nu}^{+}, u_{1,\nu}^{+})=-2\nu, \quad W(u_{0,\nu}^{-}, u_{1,\nu}^{-}) = 2\nu
\]
\end{lemma}

\begin{proof}
We make the ansatz $y(\xi)= \xi^{\frac12-\nu} (1+a(\xi))$ for $\xi>1$. Inserting this ansatz into the equation $\calH_{\nu}y=0$ yields
\[
a(\xi) = -\int_\xi^\infty \int_\xi^\zeta \eta^{-1+2\nu}\,d\eta\, U_\nu(\zeta) \zeta^{1-2\nu}(1+a(\zeta))\,d\zeta
\]
which is a Volterra equation of the form
\[
a(\xi) = \int_\xi^\infty K(\eta)(1+a(\eta))\,d\eta,\quad K(\eta) = O(\eta^{-2}) \text{\ as\ }\eta\to\infty
\]
The solution is of the form $a(\xi)=O(\xi^{-1})$ where the $O(\cdot)$ is of symbol type. This gives the solution $u_{1,\nu}^{+}(\xi)$
with the desired properties. To find $u_{0,\nu}^{+}(\xi)$, we use the reduction ansatz which yields (for some $\xi_0$ sufficiently large)
\[
u_{0,\nu}^{+}(\xi) = u_{1,\nu}^{+}(\xi)\int_{\xi_0}^\xi (u_{1,\nu}^{+}(\eta))^{-2}\,d\eta = \xi^{\frac12-\nu}(1+a(\xi)) \int_{\xi_0}^\xi \eta^{-1+2\nu}(1+a(\eta))^{-2}\,d\eta
\]
If $\nu>\frac12$, then
\[
u_{0,\nu}^{+}(\xi) = \xi^{\frac12+\nu}(1+O(\xi^{-1})) \text{\ as\ }\xi\to\infty
\]
whereas the range $0<\nu\le\frac12$ yields larger errors.
\end{proof}

We can state a very important property in analogy with the case where $\la\xi\ra V\in L^1$.
In this paper we will only need the non-resonant case.

\begin{defi}\label{def:resonance}
We say that $\calH_\nu$ has a {\em zero energy resonance} iff
\[
W_{11}:= W(u_{1,\nu}^{+}, u_{1,\nu}^{-})= 0
\]
where $u_{1,\nu}^{\pm}$ are the unique solutions from Lemma~\ref{lem:zero_en}.  This is equivalent to the existence of a nonzero solution $f$
to $\calH_\nu f=0$  so that $f(\xi)$ is asymptotic to $\xi^{\frac12 -\nu}$ as $\xi\to\infty$ and to $c\,|\xi|^{\frac12 -\nu}$ as $\xi\to-\infty$
with some constant $c\ne0$.
\end{defi}

We now perturb in energy to conclude the following. Let $\xi_0>0$ be fixed so that $u_{0,\nu}^+(\xi)>0$ for all $\xi\ge \xi_0$.

\begin{lem}\label{lemma1}
For any $\lambda\in\R$, define \begin{equation}\label{18}
u_{0,\nu}^+(\xi ,\lambda ):=u_{0,\nu}^+(\xi
)-\frac{\lambda^2}{2\nu}\int_{\xi_0}^{\xi}[u_{1,\nu}^+(\xi)u_{0,\nu}^+(\eta)-u_{1,\nu}^+(\eta)u_{0,\nu}^+(\xi)]u_{0,\nu}^+(\eta
,\lambda )\, d\eta . \end{equation}  Then
$\cH_\nu\,u_{0,\nu}^+(\cdot ,\lambda )=\lambda^2u_{0,\nu}^+(\cdot ,\lambda )$.
\end{lem}
\begin{proof}
Verify that \[(-2\nu)^{-1}[u_{1,\nu}^+(\xi)u_{0,\nu}^+(\eta)-u_{1,\nu}^+(\eta)u_{0,\nu}^+(\xi)]\] is the backward Green function of $\calH_\nu$.
\end{proof}

Next, we extend $u_{0,\nu}^+(\cdot ,\lambda )$ to a basis of solutions for $\calH_\nu f = \lambda^2 f$ for all small $\lambda>0$.

\begin{cor}\label{cor4} Let $u_{0,\nu}^+(\cdot ,\lambda )$ be defined as in \eqref{18}.
There exists  a solution $u_{1,\nu}^{+}(\cdot,\lambda)$ of $\calH_\nu f = \lambda^2 f$ with
\begin{equation}\label{eq:wu0u1}
W(u_{1,\nu}^{+}(\cdot,\lambda), u_{0,\nu}^{+}(\cdot,\lambda))=1
\end{equation}
so that for $j=0,1$  and in the range $\xi_0\le \xi \ll \lambda^{-1}$,
\begin{align}\label{33}
u_{j,\nu}^+(\xi ,\lambda )&=u_{j,\nu}^+(\xi )(1+a_{j,\nu}^+(\xi,\lambda))
\end{align}
where
\[
|\partial_\xi^\ell \partial_\lambda^k  a_{j,\nu}^+(\xi,\lambda)|\le
C_{k,\ell}\, \lambda^{2-k}\la \xi\ra^{2-\ell}
\]
for all $k,\ell\ge0$ provided $\nu>1$. In the range $0<\nu\le 1$ one has the weaker bounds
\[
|\partial_\xi^\ell \partial_\lambda^k  a_{1,\nu}^+(\xi,\lambda)|\le
\left\{\begin{array}{ll}
C_{k,\ell}\, \lambda^{2\nu-k}\la \xi\ra^{2\nu-\ell} & \nu<1 \\
C_{k,\ell}\, \lambda^{2-k}\la \xi\ra^{2-\ell}|\log(\xi\lambda)| & \nu=1
\end{array}\right.
\]
There is an analogous construction on $\xi\le0$.
\end{cor}

\begin{proof}
With $u_{0,\nu}^+(\xi ,\lambda )$ as in Lemma~\ref{lemma1}, write
\[ u_{0,\nu}^+(\xi,\lambda) = u_{0,\nu}^+(\xi)h(\xi,\lambda)\]
for all $\xi>\xi_0$.
Then
\begin{align*}
h(\xi,\lambda) &= 1 -\frac{\lambda^2}{2\nu} \int_{\xi_0}^\xi
\Big[u_{1,\nu}^+(\eta)u_{0,\nu}^+(\eta)-(u_{0,\nu}^+)^2(\eta)\frac{u_{1,\nu}^+(\xi)}{u_{0,\nu}^+(\xi)}\Big]\,
h(\eta,\lambda)\, d\eta \\
&= 1 + \lambda^2  \int_{\xi_0}^\xi K_\nu(\xi,\eta) h(\eta,\lambda)\,d\eta
\end{align*}
where $|K_\nu(\xi,\eta)|\les \eta$. Therefore, $h= 1+O(\lambda^2\xi^2)$ as claimed. For the derivatives, use
the symbol character of the $O(\cdot)$ terms from above. For $u_{1,\nu}^+$, we use the Wronskian condition:
\[ \pr_\xi\left[\frac{u_{1,\nu}^+(\xi,\lambda)}{u_{0,\nu}^+(\xi,\lambda)}\right] =
\frac{-1}{(u_{0,\nu}^+(\xi,\lambda))^2} \] and thus, with a sufficiently small constant $c$, we define
\[
u_{1,\nu}^+(\xi,\lambda) := u_{0,\nu}^+(\xi,\lambda) \int_\xi^{c\lambda^{-1}}
(u_{0,\nu}^+(\eta,\lambda))^{-2}\, d\eta .
\]
Inserting the expansion for $u_{0,\nu}^+(\eta,\lambda)$ into this
expression finishes the proof.
\end{proof}

Next, we express the Jost solutions
$f_{\pm,\nu}(\xi,\lambda)$ of $\cH_\nu$ in terms of the bases that we just constructed. Recall that the Jost solutions
 are defined as
the unique solutions to the problem
\begin{equation}
\label{eq:HnJost} \cH_\nu f_{\pm,\nu}(\cdot,\lambda) = \lambda^2
  f_{\pm,\nu}(\cdot,\lambda),\qquad f_{\pm,\nu}(\xi,\lambda)\sim
  e^{i\lambda\xi}\text{\ \  as\ \ }\xi\to\pm\infty
\end{equation}

\begin{cor}\label{cor1}
With $f_{\pm,\nu}(\cdot ,\lambda )$ being the Jost solutions of
$\cH_\nu$ with asymptotic behavior $e^{\pm i\lambda\xi}$ as
$\xi\to\pm\infty$, one has for any $\lambda\neq 0$
\begin{equation}\begin{aligned}\label{20}
f_{+,\nu}(\xi ,\lambda )&=a_{+,\nu}(\lambda )u_{0,\nu}^+(\xi ,\lambda )+b_{+,\nu}(\lambda )u_{1,\nu}^+(\xi ,\lambda )\\
f_{-,\nu}(\xi ,\lambda )&=a_{-,\nu}(\lambda )u_{0,\nu}^-(\xi ,\lambda
)+b_{-,\nu}(\lambda )u_{1,\nu}^-(\xi ,\lambda )
\end{aligned}
\end{equation} \noindent where $a_{\pm,\nu}(\lambda )= -W(f_{\pm,\nu}(\cdot
,\lambda ),u_{1,\nu}^{\pm}(\cdot ,\lambda ))$ and $b_{\pm,\nu}(\lambda
)= W(f_{\pm,\nu}(\cdot ,\lambda ),u_{0,\nu}^\pm(\cdot ,\lambda
))$.
\end{cor}

\begin{proof}
The Wronskian relations for $a_{\pm,\nu}$, $b_{\pm,\nu}$ follow
immediately from (\ref{eq:wu0u1}).
\end{proof}

To obtain an asymptotic expansion of $f_{\pm,\nu}(\cdot,\lambda)$ for large $\xi$ we introduce
\begin{equation}
\label{eq:V1_def} \cH_{0,\nu} = -\pr_\xi^2 + \big(\nu^2-\frac14\big) \xi^{-2}
\end{equation}

\begin{lem}\label{lemma4}
For any $\lambda >0$ the problem
\begin{align*}
\cH_{0,\nu} f_{0,\nu}(\cdot ,\lambda )&=\lambda^2f_{0,\nu}(\cdot ,\lambda
),\qquad  f_{0,\nu}(\xi ,\lambda )\sim e^{i\xi\lambda} \text{\ \ as\ \
}\xi\to\infty
\end{align*}
has a unique solution on $\xi >0$.  It is given by
 \begin{equation}\label{23} f_{0,\nu}(\xi ,\lambda
)=\sqrt{\frac{\pi}{2}}\, e^{i(2\nu+1)\pi
/4}\sqrt{\xi\lambda}\,H^{(+)}_\nu(\xi\lambda). \end{equation} Here
$H^{(+)}_\nu(z)=J_\nu(z)+iY_\nu(z)$ is the Hankel function.
\end{lem}

\begin{proof}
It is well-known that the ordinary differential equation
$$w''(z)+\left(\lambda^2-\frac{\nu^2-1/4}{z^2}\right)w(z)=0$$
has a fundamental system $\sqrt{z}\, J_\nu(\lambda z)$, $\sqrt{z}\,
Y_\nu(\lambda z)$ or equivalently, $\sqrt{z}\, H^{(+)}_\nu(\lambda
z)$, $\sqrt{z}\, H^{(-)}_\nu(\lambda z)$ (see~\cite{AbrSteg}).
Recall the asymptotics
\begin{equation}\label{eq:largexi}
\begin{aligned}
H^{(+)}_\nu(x)&\sim\sqrt{\frac{2}{\pi x}}\,
e^{i(x-(2\nu+1)\frac{\pi}{4})}
\qquad\mathrm{as}\; x\to +\infty\\
H^{(-)}_\nu(x)&\sim\sqrt{\frac{2}{\pi x}}\,
e^{-i(x-(2\nu+1)\frac{\pi}{4})}\qquad\mathrm{as}\; x\to +\infty.
\end{aligned}
\end{equation} Thus, (\ref{23}) is the unique solution so that $f_{0,\nu}(\xi
,\lambda )\sim e^{i\xi\lambda},$ as claimed.
\end{proof}

We now note the following:  $\{f_j(x,\lambda)\}_{j=1,2}$ is a
fundamental system of
\begin{equation}
  \label{eq:H0n_rescale}  \cH_{0,\nu} f = -f'' + \frac{\nu^2-\frac14}{x^2} f = f + \lambda
U_\nu(x,\lambda) f\text{\ \ on\ \ }x>0
\end{equation}
where $U_\nu(x,\lambda):= \lambda^{-3} U_\nu(\lambda^{-1} x)$ iff
$\{f_j(\lambda \xi,\lambda)\}_{j=1,2}$ is a fundamental system of
\[
\cH_\nu \, y(\cdot,\lambda) = \lambda^2 y(\cdot,\lambda).
\]
We shall use the bounds
\begin{equation}
  \label{eq:Wn_bds} |\partial_x^\ell \partial_\lambda^k
  U_\nu(x,\lambda)|\le C_{k,\ell}\; x^{-3-\ell} \lambda^{-k} \qquad
  \forall\; k,\ell\ge0
\end{equation}
valid for all $x\ge \lambda$.

\begin{lem}
  \label{lem:xge1} A fundamental system of \eqref{eq:H0n_rescale} on
  $x\ge \frac12$ is given by
 \[
\phi_1(x,\lambda):= \sqrt{x} H_\nu^{(+)} (x)(1+\lambda
b_1(x,\lambda)),\qquad \phi_2(x,\lambda) :=
\overline{\phi_1(x,\lambda)}
 \]
 where\footnote{This can be strengthened to $|\Im \partial_x^\ell \partial_\lambda^k b_1(x,\lambda)|\le C_{k,\ell}\,
\lambda^{-k} x^{-2-\ell}$ and $|\Re \partial_x^\ell
\partial_\lambda^k b_1(x,\lambda)|\le C_{k,\ell}\, \lambda^{-k}
x^{-3-\ell}$}
 \begin{equation}\label{eq:b1_derkl}
|\partial_x^\ell \partial_\lambda^k b_1(x,\lambda)|\le C_{k,\ell}\,
\lambda^{-k} x^{-2-\ell} \quad\forall\; k,\ell\ge0
 \end{equation}
 and all $x\ge\frac12$.
\end{lem}
\begin{proof} Let $\phi_0(x):=\sqrt{x} H_\nu^{(+)} (x)$ and observe
that $\phi_0(x)(1+\lambda b(x,\lambda))$
solves~\eqref{eq:H0n_rescale} iff \[b''(x,\lambda) \phi_0(x) +
2b'(x,\lambda)\phi_0'(x) = -U_\nu(x,\lambda)  \phi_0(x)(1+\lambda
b(x,\lambda))
\] or
\begin{equation}\label{eq:bVolterra}
b(x,\lambda) = -\int_x^\infty  \phi_0^2(y) \Big[\int_x^y
\phi_0^{-2}(u)\, du\Big] U_\nu(y,\lambda) (1+\lambda b(y,\lambda))\,
dy
\end{equation}
Define
\begin{equation}\label{eq:b0Volterra}
b_0(x,\lambda) := -\int_x^\infty  \phi_0^2(y) \Big[\int_x^y
\phi_0^{-2}(u)\, du\Big] U_\nu(y,\lambda) \, dy
\end{equation}
 Note that $\phi_0$ never vanishes and satisfies the
asymptotic expansion
\begin{equation}\label{eq:Hankel_asymp}
\phi_0(x) = c\, e^{ix} (1+ O_\R(x^{-2}) + iO_\R(x^{-1})) = e^{ix}
(1+ O_\C(x^{-1})) \text{\ \ as\ \ }x\to\infty
\end{equation}
where the $O_\R(\cdot)$ (or $O_\C(\cdot)$) terms are real-valued (or
complex-valued) and behave like symbols under differentiation. In
particular, integrating by parts shows that
\[
\sup_{y\ge x\ge\frac12} \left| \phi_0^2(y) \Big[\int_x^y \phi_0^{-2}(u)\,
du\Big] \right| \les 1
\]
which implies that $|b_0(x,\lambda)|\les x^{-2}$. The first
derivative is given by
\begin{equation}\label{eq:b0_der}
\partial_x b_0(x,\lambda) := \int_x^\infty  \phi_0^2(y)
\phi_0^{-2}(x)\, U_\nu(y,\lambda) \, dy.
\end{equation}
Since
\[
\phi_0^2(y) \phi_0^{-2}(x) = e^{2i(y-x)} (1+ O_\C(x^{-1}))
 (1+ O_\C(y^{-1}))
\]
we can integrate by parts in~\eqref{eq:b0_der} to conclude that
$|\partial_x b_0(x,\lambda)|\les x^{-3}$. For the second derivative,
\begin{equation}\label{eq:b0_der2}
\partial_x^2 b_0(x,\lambda) := \int_x^\infty \partial_x\big[ \phi_0^2(y)
\phi_0^{-2}(x)\big]\, U_\nu(y,\lambda) \, dy - U_\nu(x,\lambda)
\end{equation} we again integrate by parts
using the identity (dropping the $\C$-subscript for the remainder of
the proof)
\begin{align*}
\partial_x\big[\phi_0^2(y) \phi_0^{-2}(x)\big] &= -\partial_y\big[\phi_0^2(y)
\phi_0^{-2}(x)\big]\\
&\qquad  + e^{2i(y-x)}\Big[ O(x^{-2})
 (1+ O(y^{-1})) + O(y^{-2})
 (1+ O(x^{-1}))\Big]
\end{align*}
This yields
\begin{align}
\partial_x^2 b_0(x,\lambda) &:= \int_x^\infty \partial_x\big[ \phi_0^2(y)
\phi_0^{-2}(x)\big]\, U_\nu(y,\lambda) \, dy - U_\nu(x,\lambda)\nonumber \\
& = \int_x^\infty \phi_0^2(y) \phi_0^{-2}(x)\, \partial_y
U_\nu(y,\lambda) \, dy + O (x^{-4}).\label{eq:Wn_int2}
\end{align}
Integrating by parts in the integral on line~\eqref{eq:Wn_int2} we
conclude that $|\partial_x^2 b_0(x,\lambda)|\les x^{-4}$. Continuing
 in this fashion one proves that for all $\ell\ge0$,
\[
|\partial_x^\ell b_0(x,\lambda)|\le C_\ell\, x^{-2-\ell}
\]
In view of \eqref{eq:Wn_bds}, the $\lambda$-derivatives are treated
in exactly the same way and we thus obtain the estimates
\begin{equation}\label{eq:b0_derkl}
|\partial_x^\ell\partial_\lambda^k\,  b_0(x,\lambda)|\le
C_{\ell,k}\, x^{-2-\ell}\lambda^{-k}
\end{equation}
for all $k,\ell\ge0$. These estimates transfer via
\eqref{eq:bVolterra} to $b(x,\lambda)$ because $\lambda$ is small.
Indeed, first note that~\eqref{eq:bVolterra} has a solution via a
contraction, say. Second, repeating the same arguments that lead
to~\eqref{eq:b0_derkl} but with~\eqref{eq:b1_derkl} as bootstrap
assumptions shows that we can get~\eqref{eq:b1_derkl} back with the
same constants (provided those are sufficiently large); the point is
as follows: estimating $\partial_x^\ell b(x,\lambda)$ requires at
most $\partial_y^{\ell-1} b(y,\lambda)$ inside the integral on the
right-hand side of~\eqref{eq:bVolterra} (see \eqref{eq:Wn_int2} for
the case $\ell=2$). While $\partial_\lambda^k b(x,\lambda)$ with $k$
fixed can appear on both sides of~\eqref{eq:bVolterra}, note that
then we can use $\lambda$ small (with a smallness that does not
depend on~$k$) to solve for that derivative.
\end{proof}

Next, we describe a basis of solutions for $\lambda\ll x\le \frac12$.

\begin{lem}
\label{lem:xle1} A fundamental system of
\eqref{eq:H0n_rescale} on
  $\lambda\ll x\le\frac12$ is given by
 \[
\psi_1(x,\lambda):= \sqrt{x} J_\nu(x)(1+\lambda
c_1(x,\lambda)),\qquad \psi_2(x,\lambda) := \sqrt{x}\,
Y_\nu(x)(1+\lambda c_2(x,\lambda))
 \]
where for $j=1,2$, $c_j(x,\lambda)$ are real-valued and satisfy the
bounds
 \begin{equation}\label{eq:cj_derkl}
|\partial_x^\ell \partial_\lambda^k c_j(x,\lambda)|\le C_{k,\ell}\,
\lambda^{-k} x^{-1-\ell} \quad\forall\; k,\ell\ge0
 \end{equation}
 and all $\lambda\ll x\le\frac12$.
\end{lem}
\begin{proof} As in the previous lemma, and since $Y_\nu(x)<0$ for
all $0<x\le \frac12$ provided $\nu\ge0$,
\begin{equation}\label{eq:c2Y}
c_2(x,\lambda) = \int_x^1 yY_\nu^2(y) \int_x^y u^{-1}Y_\nu^{-2}(u)\,
du\;U_\nu(y,\lambda)(1+\lambda c_2(y,\lambda))\, dy
\end{equation}
Recall the asymptotic behavior, as $x\to0+$ and with real constants
$\alpha_{1,\nu}, \alpha_{2,\nu}$,
\begin{equation}
  \label{eq:JnuYnu} J_\nu(x) = \alpha_{1,\nu}\, x^\nu(1+x\omega_1(x)),\qquad Y_\nu(x)
= \alpha_{2,\nu}\, x^{-\nu}(1+x\omega_2(x))
\end{equation}
where $\omega_j$ behave like symbols under differentiation: for all
$\ell\ge0$ and $j=1,2$,
\[
|\omega_j^{(\ell)}(x)|\le C_\ell\, x^{-\ell}, \qquad\; 0<x<\frac12
\]
 First,
let
\[
c_{2,0}(x,\lambda) := \int_x^{\frac12} yY_\nu^2(y) \int_x^y
u^{-1}Y_\nu^{-2}(u)\, du\;U_\nu(y,\lambda)\, dy
\]
Then for all $x>\lambda$,
\begin{align*}
  |c_{2,0}(x,\lambda)| &\les \int_x^{\frac12} y^{1-2\nu} \int_x^y
u^{-1+2\nu}\, du\;  |U_\nu(y,\lambda)|\, dy \les \int_x^\infty y^{-2}\,
dy\les
x^{-1} \\
|\pr_x c_{2,0}(x,\lambda)| &\les \int_x^{\frac12} yY_\nu^2(y)
x^{-1}Y_\nu^{-2}(x) |U_\nu(y,\lambda)|\, dy \les \int_x^{\frac12} y^{1-2\nu}
x^{-1+2\nu} y^{-3}\, dy \les x^{-2}
\end{align*}
Inductively, it now follows that
\begin{equation}\label{eq:c20_der}
|\pr_x^\ell \partial_\lambda^k c_{2,0}(x,\lambda)| \le C_{k,\ell}\,
\lambda^{-k} x^{-1-\ell} \quad\forall\; k,\ell\ge0
\end{equation}
By a fixed-point argument, \eqref{eq:c2Y} has a solution
$c_2(x,\lambda)$ for small $\lambda$ on $\lambda\le x\le \frac12$ which
satisfies~\eqref{eq:cj_derkl} for all $k\ge0$ and $\ell=0$. The same
arguments that lead to~\eqref{eq:c20_der} now
yield~\eqref{eq:cj_derkl} for all $\ell>0$ and $j=2$, settling the
case of~$\psi_2$.

\noindent The solution $\psi_1(x,\lambda)$ is given by, with a
suitable constant $\gamma_\nu\ne0$,
\begin{equation}\label{eq:psi1}\begin{split}
 \psi_1(x,\lambda) &:= \gamma_\nu^{-1} \psi_2(x,\lambda) \int_{0}^x
\psi_2^{-2}(y,\lambda)\, dy  \\
&= \sqrt{x}\, Y_\nu(x) (1+\lambda c_2(x,\lambda)) \int_{0}^x
y^{-1}Y_\nu(y)^{-2} (1+\lambda c_2(y,\lambda))^{-2}\, dy
\end{split}
\end{equation}
for some sufficiently large constant $A$ which insures that
\[\lambda c_2\ll 1 \text{\ \ on\ \ } A\lambda\le x\le \frac12\] Moreover, we set
$c_2(x,\lambda):=0$ for all $0\le x\le A\lambda$. Due to this fact,
$\psi_1$ as defined in~\eqref{eq:psi1} solves~\eqref{eq:H0n_rescale}
only on the interval $A\lambda\le x\le \frac12$, which however is
sufficient for our purposes. The  constant $\gamma_\nu\ne0$ is
defined via the relation
\[
\gamma_\nu\, \sqrt{x} J_\nu(x) =  \sqrt{x}\, Y_\nu(x) \int_{0}^x
y^{-1}Y_\nu(y)^{-2} \, dy
\]Hence,
we see that
\begin{align*}
  & \gamma_\nu^{-1}\, \sqrt{x}\, Y_\nu(x) (1+\lambda c_2(x,\lambda)) \int_{0}^x
y^{-1}Y_\nu(y)^{-2} (1+\lambda c_2(y,\lambda))^{-2}\, dy \\
 & =  \sqrt{x} J_\nu(x) (1+\lambda c_2(x,\lambda))\Big[ 1
 + \lambda Y_\nu(x) J_\nu^{-1}(x)\int_0^x y^{-1}Y_\nu(y)^{-2} O(c_2(y,\lambda))\,
 dy\Big]\\
 & =  \sqrt{x} J_\nu(x) (1+\lambda c_2(x,\lambda))\Big[ 1
 + \lambda O(x^{-2\nu} x^{-1+2\nu})\Big] \\
 & =: \sqrt{x} J_\nu(x) (1+\lambda c_1(x,\lambda))
\end{align*}
with $c_1$ inheriting the bounds \eqref{eq:cj_derkl} from $c_2$.
\end{proof}

In what follows, $\beta_\nu:= \sqrt{\frac{\pi}{2}}\, e^{i(2\nu+1)\pi
/4}$. We can now describe the Jost solutions of $\cH_\nu$ in the
region $1\ll\xi\le \lambda^{-1}$, which is needed for the matching
described in Corollary~\ref{cor1}.

\begin{cor}
  \label{cor:regionI} For $\lambda\ne0$, let $f_{+,\nu}(\xi,\lambda)$ be the Jost
  solution satisfying~\eqref{eq:HnJost}.  Then for all $1\ll\xi\le
  \lambda^{-1}$ there is the representation
  \begin{equation}
    \label{eq:regionIJost} \begin{aligned}
      f_{+,\nu}(\xi,\lambda) &= \beta_\nu\,\sqrt{\lambda\xi}\big[
      J_\nu(\lambda\xi) (1+O(\lambda))(1+O(\xi^{-1})) +
      Y_\nu(\lambda\xi) O(\lambda)(1+O(\xi^{-1}))\big] \\
      &\quad +i\beta_\nu\,\sqrt{\lambda\xi}\big[
      Y_\nu(\lambda\xi) (1+O(\lambda))(1+O(\xi^{-1})) +
      J_\nu(\lambda\xi) O(\lambda)(1+O(\xi^{-1}))\big]
    \end{aligned}
  \end{equation}
where each $O(\lambda)$ and $O(\xi^{-1})$ is real-valued and behaves
like a symbol under differentiation:
\begin{equation}\label{eq:Odiff}
|\partial_\xi^\ell \partial_\lambda^k O(\xi^{-1})|\le C_{k,\ell}\,
\lambda^{-k} \xi^{-1-\ell}, \quad  |\partial_\xi^\ell
\partial_\lambda^k O(\lambda)|\le C_{k,\ell}\, \lambda^{1-k}
\xi^{-\ell} \quad\forall\; k,\ell\ge0
\end{equation}
in the range $0<\lambda\ll1$, $1\ll \xi\le \lambda^{-1}$.
\end{cor}
\begin{proof}
  In the rescaled picture, i.e., with $x=\lambda\xi$ the
  corresponding representation is given by
  \[
\begin{aligned}
       & \beta_\nu\,\sqrt{x}\big[
      J_\nu(x) (1+O(\lambda))(1+\lambda c_1(x,\lambda)) +
      Y_\nu(x) O(\lambda)(1+\lambda c_2(x,\lambda))\big] \\
      & +i\beta_\nu\,\sqrt{x}\big[
      Y_\nu(x) (1+O(\lambda))(1+\lambda c_2(x,\lambda)) +
      J_\nu(x) O(\lambda)(1+\lambda c_1(x,\lambda))\big]
    \end{aligned}
  \]
  as can be seen by matching the solutions of Lemma~\ref{lem:xle1} to
  those of Lemma~\ref{lem:xge1} at $x=1$. The $O(\lambda)$ have the
  claimed symbol behavior due to~\eqref{eq:b1_derkl}. Furthermore,
  \[
\lambda c_j(\lambda \xi,\lambda) = O(\xi^{-1})
  \]
  behaves under differentiation as claimed, see~\eqref{eq:cj_derkl},
  and we are done.
\end{proof}

The Wronskians appearing in Corollary~\ref{cor1} will be evaluated
at $\xi=\lambda^{-1+\eps}$ where $\eps>0$ behaves like
$\frac{1}{4\nu}$. As a preliminary step, we note the following for
$0<\lambda\ll1$:

\begin{cor}
  \label{cor:match_region} For $\lambda\ne0$, let $f_{+,\nu}(\xi,\lambda)$ be the Jost
  solution of $\cH_\nu$.  Then for sufficiently small $\eps>0$ and all $\lambda^{-1+\eps}\le \xi\le
  \lambda^{-1}$ there is the representation
  \begin{equation}
    \label{eq:regionIJost'} \begin{aligned}
      f_{+,\nu}(\xi,\lambda) &= \beta_\nu\,\sqrt{\lambda\xi}\,
      J_\nu(\lambda\xi) \big[1+O(\lambda^\eps)\big]
      +i\beta_\nu\,\sqrt{\lambda\xi}\,\,
      Y_\nu(\lambda\xi) \big[1+ O(\lambda^{1-\eps})\big]
    \end{aligned}
  \end{equation}
where $O(\lambda^\eps)$ and $O(\lambda^{1-\eps})$ are real-valued
and behave like symbols under differentiation:
\[
|\partial_\xi^\ell \partial_\lambda^k O(\lambda^\eps)|\le
C_{k,\ell}\, \lambda^{\eps-k} \xi^{-\ell},\qquad |\partial_\xi^\ell
\partial_\lambda^k O(\lambda^{1-\eps})|\le C_{k,\ell}\, \lambda^{1-\eps-k}
\xi^{-\ell}  \quad\forall\; k,\ell\ge0
\]
in the range $0<\lambda\ll1$, $\lambda^{-1+\eps}\le \xi\le
\lambda^{-1}$.
\end{cor}
\begin{proof} Simply note that, for $\eps\le \frac{1}{4\nu}$,
\[
\lambda \frac{Y_\nu(\lambda\xi)}{ J_\nu(\lambda\xi)} = \lambda O(
(\lambda\xi)^{-2\nu}) = O(\lambda^\eps)
\]
in the specified range. Hence, \eqref{eq:regionIJost'} follows from
\eqref{eq:regionIJost}. The behavior under differentiation is also
clear from~\eqref{eq:Odiff}.
\end{proof}

We now can compute the coefficients $a_{\pm,\nu}(\lambda),
b_{\pm,\nu}(\lambda)$ and the Wronskian $W_\nu(\lambda)$ for small
$\lambda$, see Corollary~\ref{cor1}. Recall that
$O(\lambda^\sigma)$ behaves like a symbol  under differentiation if
\[ \pr_\lambda^\ell O(\lambda^\sigma) = O(  \lambda^{\sigma-\ell}) \text{\ \ as\ \ }\lambda\to0+\]
for all $\ell\ge0$.

\begin{prop}
  \label{prop:wronski} Let $\beta_\nu$ be as above. With nonzero real constants
  $\alpha_{0,\nu}^+,\beta_{0,\nu}^+$, and some
  sufficiently small $\eps>0$,
  \begin{equation}
    \label{eq:abW}
    \begin{aligned}
      a_{+,\nu}(\lambda) &= \lambda^{\frac12+\nu} \beta_\nu(\alpha_{0,\nu}^+ +
      O(\lambda^\eps) + iO(\lambda^{(1-2\nu)\eps})) \\
 b_{+,\nu}(\lambda) &= i\lambda^{\frac12-\nu} \beta_\nu(\beta_{0,\nu}^+ +
      O(\lambda^\eps) + iO(\lambda^{(1+2\nu)\eps}))
    \end{aligned}
  \end{equation} as $\lambda\to0+$
  with real-valued $O(\cdot)$ which behave like symbols under
  differentiation in~$\lambda$. The asymptotics as $\lambda\to0-$
  follows from that as $\lambda\to0+$ via the relations
  $a_{+,\nu}(-\lambda)=\overline{a_{+,\nu}(\lambda)}$,
  $b_{+,\nu}(-\lambda)=\overline{b_{+,\nu}(\lambda)}$. Analogous expressions hold for $a_{-,\nu}$
and $b_{-,\nu}$.
\end{prop}

\begin{proof} Evaluating at $\xi=\lambda^{-1+\eps}$, using
\eqref{eq:JnuYnu} and Lemma~\ref{lem:zero_en},
  \begin{align*}
    W(f_{+,\nu}(\cdot,\lambda), u_{1,\nu}^+(\cdot,\lambda)) &= \beta_\nu W\big( \sqrt{\lambda\xi}\,
      J_\nu(\lambda\xi) \big[1+O(\lambda^\eps)\big],
      u_{1,\nu}^+(\xi)(1+O(\lambda^2\xi^2))\big) \\
      &\quad + i\beta_\nu W\big( \sqrt{\lambda\xi}\,
      Y_\nu(\lambda\xi) \big[1+O(\lambda^\eps)\big],
      u_{1,\nu}^+(\xi)(1+O(\lambda^2\xi^2))\big) \\
 &= \beta_\nu \alpha_{1,\nu}\, \lambda^{\frac12+\nu} W\big( \xi^{\frac12+\nu} \big[1+O(\lambda^\eps)\big],
      \xi^{\frac12-\nu}(1+O(\lambda^\eps))\big) \\
      &\quad + i\beta_\nu \alpha_{2,\nu}  \lambda^{\frac12-\nu} W\big( \xi^{\frac12-\nu} \big[1+O(\lambda^\eps)\big],
      \xi^{\frac12-\nu} (1+O(\lambda^\eps))\big) \\
      &= \beta_\nu\, \lambda^{\frac12+\nu} [\alpha_{1,\nu}+O(\lambda^\eps) + i
      O((\xi \lambda)^{-2\nu} \lambda^\eps)]
  \end{align*}
  with some constant $\tilde\beta_\nu\ne0$. Next,
  \begin{align*}
    W(f_{+,\nu}(\cdot,\lambda), u_{0,\nu}^+(\cdot,\lambda)) &= \beta_\nu W\big( \sqrt{\lambda\xi}\,
      J_\nu(\lambda\xi) \big[1+O(\lambda^\eps)\big],
      u_{0,\nu}^+(\xi)(1+O(\lambda^2\xi^2))\big) \\
      &\quad + i\beta_\nu W\big( \sqrt{\lambda\xi}\,
      Y_\nu(\lambda\xi) \big[1+O(\lambda^\eps)\big],
      u_{0,\nu}^+(\xi)(1+O(\lambda^2\xi^2))\big) \\
 &= \beta_\nu \alpha_{1,\nu} \lambda^{\frac12+\nu} W\big( \xi^{\frac12+\nu} \big[1+O(\lambda^\eps)\big],
      \xi^{\frac12+\nu}(1+O(\lambda^\eps))\big) \\
      &\quad + i\beta_\nu \alpha_{2,\nu}\lambda^{\frac12-\nu} W\big( \xi^{\frac12-\nu} \big[1+O(\lambda^\eps)\big],
      \xi^{\frac12+\nu} (1+O(\lambda^\eps))\big) \\
      &= i\beta_\nu \, \lambda^{\frac12-\nu} [\alpha_{2,\nu} +O(\lambda^\eps) +
      i
      O((\xi \lambda)^{2\nu}\lambda^{\eps})]
  \end{align*}
The proposition now follows by combining these calculations with
Corollary~\ref{cor1}.
\end{proof}

We can now describe the Wronskian $W_\nu(\lambda)=W(f_{-,\nu}(\cdot,\lambda),f_{+,\nu}(\cdot,\lambda))$ in the non-resonant case:

\begin{cor}\label{cor:Wnu_nonres}
If $\calH_\nu$ is nonresonant in the sense of Definition~\ref{def:resonance}, then with some small $\eps>0$ depending on~$\nu$,
\begin{equation}\label{eq:Wnu_asymp}
\begin{aligned}
W_\nu(\lambda) &=  b_{+,\nu}(\lambda) b_{-,\nu}(\lambda) (W_{11}+O(\lambda^\eps)+ iO(\lambda^{(1-2\nu)\eps})) \\
&=i e^{i\nu\pi}\, \lambda^{1-2\nu} (W_{0,\nu} + O_{\C}(\lambda^\eps)) \text{\ \ as\ \ }\lambda\to0+
\end{aligned}\end{equation}
Here $W_{0,\nu}$  is a  nonzero real constant and $O_\C(\lambda^\eps)$ is complex valued and  of symbol type.
For $\lambda<0$, one has $W_\nu(-\lambda)=\overline{W_\nu(\lambda)}$.
\end{cor}
\begin{proof}
In view of \eqref{20},
\begin{equation}\nonumber
\begin{aligned}
W_\nu(\lambda) &:= W(f_{-,\nu}(\cdot,\lambda), f_{+,\nu}(\cdot,\lambda))\\
& = a_{-,\nu}(\lambda )a_{+,\nu}(\lambda ) W(u_{0,\nu}^-(\cdot ,\lambda
), u_{0,\nu}^+(\cdot ,\lambda
)) + a_{-,\nu}(\lambda ) b_{+,\nu}(\lambda ) W(u_{0,\nu}^-(\cdot ,\lambda
), u_{1,\nu}^+(\cdot ,\lambda )) \\
&+ b_{-,\nu}(\lambda )a_{+,\nu}(\lambda ) W(u_{1,\nu}^-(\cdot ,\lambda
), u_{0,\nu}^+(\cdot ,\lambda
)) +    b_{-,\nu}(\lambda ) b_{+,\nu}(\lambda ) W(u_{1,\nu}^-(\cdot ,\lambda
), u_{1,\nu}^+(\cdot ,\lambda
))
\end{aligned}\end{equation}
for all $\lambda\ne0$. From \eqref{33},
\[
W(u_{1,\nu}^-(\cdot ,\lambda
), u_{1,\nu}^+(\cdot ,\lambda
)) = W_{11} + O(\lambda^\alpha) \text{\ as\ }\lambda\to0+
\]
where $0<\alpha= \min(1,2\nu)$ (with a logarithmic loss at $\nu=1$). By our nonresonant assumption, $W_{11}\ne0$.
Using the asymptotic expansions of Proposition~\ref{prop:wronski}
as well as setting
\[
W(u_{j,\nu}^-(\cdot ,\lambda
), u_{k,\nu}^+(\cdot ,\lambda )) = O(1)
\]
for all $j+k<1$, finishes the proof.
\end{proof}


\section{The scattering theory of $\calH_{d,n}$, $d+n>1$}
\label{sec:scatter2}

As explained in Section~\ref{sec:setup}, we reduce the Laplacean on $\calM$ to the Schr\"odinger operator
\[
\calH_{d,n} = -\partial_\xi^2 + \big(2\mu_n^2 + d(d-2)/4)\la\xi\ra^{-2} + O(\la\xi\ra^{-3})
\]
where the $O(\cdot)$ is of symbol type. This means that $\calH_{d,n}=\calH_\nu$ in the sense of Definition~\ref{def:Hnu} with
\[
\nu = \sqrt{2\mu_n^2 + (d-1)^2/4}
\]
Note that $\nu>0$ unless $d=1, n=0$ which is not allowed here (this case was considered in Part~I, see~\cite{SSS}).
To be able to apply the results of Section~\ref{sec:scatternu}, we need to verify the following:

\begin{lemma}\label{lem:legit}
For any $d+n>1$ the operator
$\calH_{d,n}$ does not
have a {\em zero energy resonance} in the sense of Definition~\ref{def:resonance}.
\end{lemma}
\begin{proof}
If there were a solution $u$ of $\calH_{d,n}u=0$ with
the property that $|u(\xi)|\les |\xi|^{\frac12-\nu}$ as $\xi\to\pm\infty$, then lifting
this to $\calM$ would yield a harmonic function decaying like
$|\xi|^{1/2-2\nu-d/2}$ at both ends. But since $d\ge1,$ this would imply the existence of a
nonzero harmonic function on~$\calM$ that vanishes at both ends (for all $\nu>0$). However, by the maximum principle such a
harmonic function would need to vanish identically. This contradiction rules out a zero energy resonance of $\calH_{d,n}$.
\end{proof}

Thus, Proposition~\ref{prop:wronski} and Corollary~\ref{cor:Wnu_nonres} apply to all $\calH_{d,n}$ with $d+n>1$.

\section{The oscillatory  integral estimates for $d+n>1$}
\label{sec:oscill}

In this section we associate $\calH_{d,n}$ with $\calH_\nu$ as in the previous section.
In fact, the estimates of this section do not use any other information about $\calH_{d,n}$ than
that furnished by Section~\ref{sec:scatternu}.

We begin with a corollary to Proposition~\ref{prop:wronski}. As
already mentioned in Section~\ref{sec:setup}, the importance of this
corollary lies with the fact that the spectral resolution
$E_\nu(d\lambda)(\xi,\xi')$ of $\cH_\nu$ satisfies, for $\xi>\xi'$,
\[
E_\nu(d\lambda^2)(\xi,\xi') = 2\lambda \Im \Big[
\frac{f_{+,\nu}(\xi,\lambda)
    f_{-,\nu}(\xi',\lambda)}{W_\nu(\lambda)} \Big]\,d\lambda
\]
as an identity of  Schwartz kernels. In this section we
prove~\eqref{eq:crux}. 
We break this proof up
into a small and a large energy piece. We also need to
distinguish the oscillatory regime from the exponential regime in
the Jost solutions~$f_{\pm,\nu}(\xi,\lambda)$ (the transition happens
at $|\lambda\xi|=1$). This section will freely use the notations of Section~\ref{sec:scatternu}.

\begin{cor}
  \label{cor:resolv} For $0<\lambda\ll1$ and any $\xi,\xi'\in\R$,
  \begin{equation}
    \label{eq:spec_dens}\begin{split}
     & \Im \Big[   \frac{f_{+,\nu}(\xi,\lambda)
    f_{-,\nu}(\xi',\lambda)}{W_\nu(\lambda)} \Big] = O(\lambda^{2\nu} )
    u_{0,\nu}^+ (\xi,\lambda)u_{1,\nu}^-(\xi',\lambda) + O(\lambda^{2\nu}) u_{1,\nu}^+(\xi,\lambda)
    u_{0,\nu}^-(\xi',\lambda)
    \\
    & \qquad\qquad + O(\lambda^{2\nu})
    u_{0,\nu}^+(\xi,\lambda) u_{0,\nu}^-(\xi',\lambda) + O(\lambda^{2\nu}) u_{1,\nu}^+(\xi,\lambda) u_{1,\nu}^-(\xi',\lambda)
    \end{split}
  \end{equation}
  where the $O(\cdot)$ are real-valued and behave like symbols under
  differentiation in~$\lambda$.
\end{cor}
\begin{proof}
By Corollary~\ref{cor1} one has
\begin{align}
& \Im \Big[   \frac{f_{+,\nu}(\xi,\lambda)
    f_{-,\nu}(\xi',\lambda)}{W_\nu(\lambda)} \Big]  \nonumber \\
    & = \Im \Big[\frac{ (a_{+,\nu}(\lambda) u_{0,\nu}^+(\xi,\lambda) + b_{+,\nu}(\lambda) u_{1,\nu}^+(\xi,\lambda))
     (a_{-,\nu}(\lambda) u_{0,\nu}^-(\xi',\lambda) + b_{-,\nu}(\lambda) u_{1,\nu}^-(\xi',\lambda))}{W_\nu(\lambda)}
     \Big] \nonumber \\
     & = \Im
     \Big[\frac{a_{+,\nu} a_{-,\nu}(\lambda)}{W_\nu(\lambda)}\Big] \, u_{0,\nu}^+(\xi,\lambda)u_{0,\nu}^-(\xi',\lambda)
      + \Im
     \Big[\frac{b_{+,\nu}b_{-,\nu}(\lambda)}{W_\nu(\lambda)}\Big]\, u_{1,\nu}^+(\xi,\lambda)
     u_{1,\nu}^-(\xi',\lambda) \nonumber \\
     &\quad + \Im\Big[
     \frac{a_{+,\nu}b_{-,\nu}(\lambda)}{W_\nu(\lambda)} \Big]\,
     u_{0,\nu}^+(\xi,\lambda) u_{1,\nu}^-(\xi',\lambda) + \Im\Big[
     \frac{a_{-,\nu}b_{+,\nu}(\lambda)}{W_\nu(\lambda)} \Big]
     u_{0,\nu}^-(\xi,\lambda) u_{1,\nu}^+(\xi',\lambda).  \label{eq:drei_Im}
\end{align}
One first verifies from \eqref{eq:Wnu_asymp} that
\[
W_\nu(\lambda) = b_{+,\nu}(\lambda) b_{-,\nu}(\lambda) \big[1+ O(\lambda^\eps) + i\lambda^{2\nu} (\tau_\nu + O(\lambda^\eps) + iO(\lambda^{(1-2\nu)\eps})) \big]
\]
with real-valued $O(\cdot)$ terms and some real constant $\tau_\nu$.
The four imaginary parts in~\eqref{eq:drei_Im} are now computed systematically from Proposition~\ref{prop:wronski}
and this expression. For example,
\[
\Im
     \Big[\frac{b_{+,\nu}b_{-,\nu}(\lambda)}{W_\nu(\lambda)}\Big] = \Im \big[1+O(\lambda^\eps) + i\lambda^{2\nu} (\tau_\nu + O(\lambda^\eps
+ iO(\lambda^{(1-2\nu)\eps}))\big]^{-1}
= O(\lambda^{2\nu})
\]
and
\[
\Im\Big[
     \frac{a_{+,\nu}b_{-,\nu}(\lambda)}{W_\nu(\lambda)} \Big] = \Im \frac{\lambda^{\frac12+\nu} (\alpha_{0,\nu}^+
+ O(\lambda^{\eps}) +i O(\lambda^{(1-2\nu)\eps}))}{i\lambda^{\frac12-\nu} (\beta_{0,\nu}^+ + O(\lambda^\eps) + i O(\lambda^{(1+2\nu)\eps}))} = O(\lambda^{2\nu})
\]
as claimed. We leave the other two imaginary parts to the reader.
\end{proof}

We now
proceed to our first oscillatory integral estimate. Let $\chi$ be a
smooth cut-off function to small energies, i.e., $\chi (\lambda )=1$
for small $|\lambda |$ and $\chi$ vanishes outside a small interval
around zero.  In addition, we introduce the smooth cut-off functions
$\chi_{[|\xi\lambda |<1]}$ and $\chi_{[|\xi\lambda |>1]}$ which form
a partition of unity adapted to these intervals. For the remainder
of the paper, constants implicit in the $\les$ notation of course do
 not depend on~$t$.

\begin{lem}\label{lem:osc_int1}
For all $t>0$ and any $0\le\sigma\le\nu-\frac{d-1}{2}$,
\begin{align}\label{eq:disp_1} \sup_{\xi,\xi
'}\bigg|\int_0^{\infty}e^{it\lambda^2}\lambda
\chi(\lambda;\xi,\xi')(\langle\xi\rangle\langle\xi
'\rangle)^{-\frac{d}{2}-\sigma}\;\Im\left[\frac{f_{+,\nu}(\xi,\lambda
)f_{-,\nu}(\xi
',\lambda )}{W_\nu(\lambda )}\right]\, d\lambda \bigg| &\les t^{-\frac{d+1}{2}-\sigma} \\
\label{eq:disp_1wave} \sup_{\xi,\xi
'}\bigg|\int_0^{\infty}e^{\pm it\lambda}\lambda
\chi(\lambda;\xi,\xi')(\langle\xi\rangle\langle\xi
'\rangle)^{-\frac{d}{2}-\sigma}\;\Im\left[\frac{f_{+,\nu}(\xi,\lambda
)f_{-,\nu}(\xi
',\lambda )}{W_\nu(\lambda )}\right]\, d\lambda \bigg| &\les t^{-\frac{d}{2}-\sigma}
\end{align}
where $\chi(\lambda;\xi,\xi'):=\chi (\lambda )\chi_{[|\xi\lambda
|<1,|\xi '\lambda |<1]}$.
\end{lem}

\begin{proof}
We now write Corollary~\ref{cor:resolv} schematically in the form
\[
\Im \Big[   \frac{f_{+,\nu}(\xi,\lambda)
    f_{-,\nu}(\xi',\lambda)}{W_\nu(\lambda)} \Big] = O(\lambda^{2\nu}) O((\la \xi\ra\la \xi'\ra)^{\frac12+\nu})
\]
where the second $O(\cdot)$ term is obtained from Corollary~\ref{cor4}. Under differentiation in~$\lambda$ the right-hand
side behaves like a symbol.
Thus, \eqref{eq:disp_1} reduces to the following stationary phase bound
\[
\Big| \int_0^\infty e^{it\lambda^2} O(\lambda^{2\nu+1})
  \frac{\chi(\lambda;\xi,\xi')}{(\langle\xi\rangle\langle\xi
'\rangle)^{\frac{d}{2}+\sigma}}  O((\la \xi\ra\la \xi'\ra)^{\frac12+\nu}) \, d\lambda\Big|  \les t^{-\frac{d+1}{2}-\sigma}.
\]
Observe that
\[
O((\la \xi\ra\la \xi'\ra)^{\nu-\sigma-\frac{d-1}{2}}) \lambda^{2\nu-2\sigma-(d-1)} = O(1)
\]
on the support of the integrand (with symbol behavior under differentiation in $\lambda$). Hence, we conclude that
it suffices to prove
\[
 \Big|\int_0^\infty e^{it\lambda^2} \lambda^{d+2\sigma}
 \omega(\lambda;\xi,\xi')\,d\lambda \Big|\les t^{-\frac{d+1}{2}-\sigma}
\]
where for all $N\ge1$
\[
\sup_{\xi,\xi'}|\partial_\lambda^N \omega(\lambda;\xi,\xi') |\le
C_\nu\, \lambda^{-N}
\]
uniformly in $\xi,\xi'$.  However, this is a standard estimate and~\eqref{eq:disp_1} follows.

The bound \eqref{eq:disp_1wave}  is a consequence of the oscillatory integral estimate
\[
 \Big|\int_0^\infty e^{\pm it\lambda} \lambda^{d+2\sigma}
 \omega(\lambda;\xi,\xi')\,d\lambda \Big|\les \min(1,t^{-d-1-2\sigma})
\]
and the lemma follows.
\end{proof}

Next, we consider the case $|\xi\lambda |>1$ and $|\xi '\lambda
|>1$.  With the convention that $f_{\pm}(\xi,-\lambda
)=\overline{f_{\pm}(\xi,\lambda )}$ we can remove the imaginary part
from the resolvent and integrate $\lambda$ over the whole axis.
Indeed, if $a(\lambda;\xi,\xi')$ is an even function in~$\lambda$,
then
\begin{align*}
  &\int_0^\infty e^{it\lambda^2}\, \lambda \chi(\lambda)
  a(\lambda;\xi,\xi') \Im \Big[ \frac{f_{+,\nu}(\xi,\lambda)
  f_{-,\nu}(\xi',\lambda)}{W_\nu(\lambda)}\Big]\, d\lambda \\
  &= \int_0^\infty e^{it\lambda^2}\, \lambda \chi(\lambda)
  a(\lambda;\xi,\xi')  \Big[ \frac{f_{+,\nu}(\xi,\lambda)
  f_{-,\nu}(\xi',\lambda)}{W_\nu(\lambda)} - \frac{f_{+,\nu}(\xi,-\lambda)
  f_{-,\nu}(\xi',-\lambda)}{W_\nu(-\lambda)} \Big] \, d\lambda \\
  &= \int_{-\infty}^\infty e^{it\lambda^2}\, \lambda \chi(\lambda)
  a(\lambda;\xi,\xi')  \frac{f_{+,\nu}(\xi,\lambda)
  f_{-,\nu}(\xi',\lambda)}{W_\nu(\lambda)}\, d\lambda
\end{align*}

We shall follow this convention henceforth. To estimate the
oscillatory integrals, we shall repeatedly use the following version
of stationary phase, see Lemma~2 in \cite{Sch} for the proof.

\begin{lem}\label{lemma12'}
Let $\phi (0)=\phi '(0)=0$ and $1\leq\phi ''\leq C$.  Then
\begin{equation}\label{46} \bigg|\int_{-\infty}^{\infty}e^{it\phi (\lambda)}a(\lambda)\,
d\lambda\bigg|\les\delta^2\left\{\int\frac{|a(\lambda)|}{\delta^2+|\lambda|^2}\,
d\lambda+\int_{|\lambda|>\delta}\frac{|a'(\lambda)|}{|\lambda|}\,
d\lambda\right\}
\end{equation} where $\delta =t^{-1/2}$.
\end{lem}

Before proceeding, let us note that
\[
\sup_{\xi,\xi'}\bigg|\int_{-\infty}^{\infty}e^{it\lambda^2}\lambda\chi (\lambda
)\chi_{[|\xi\lambda |>1,|\xi '\lambda
|>1]}(\langle\xi\rangle\langle\xi
'\rangle)^{-\frac{d}{2}-\sigma}\frac{f_{+,\nu}(\xi,\lambda )f_{-,\nu}(\xi ',\lambda
)}{W_\nu(\lambda )}\, d\lambda \bigg| \les 1
\]
due to the fact that
\[
\sup_{\xi,\xi'} |\lambda| \Big| \frac{f_{+,\nu}(\xi,\lambda )f_{-,\nu}(\xi ',\lambda
)}{W_\nu(\lambda )} \Big|\les 1
\]
see Section~\ref{sec:scatternu}. Hence, in all small energy oscillatory integrals
it suffices to consider $t>1$. The same comment applies of course the the wave equation.

\noindent Using Lemma~\ref{lemma12'} we can prove the following:
\begin{lemma}\label{lemma13}
For all $t\ge1$ and $0\le\sigma\le\nu-\frac{d-1}{2}$,
\begin{align}\label{eq:disp_2} \sup_{\xi >0>\xi
'}\bigg|\int_{-\infty}^{\infty}e^{it\lambda^2}\lambda\chi (\lambda
)\chi_{[|\xi\lambda |>1,|\xi '\lambda
|>1]}(\langle\xi\rangle\langle\xi
'\rangle)^{-\frac{d}{2}-\sigma}\frac{f_{+,\nu}(\xi,\lambda )f_{-,\nu}(\xi ',\lambda
)}{W_\nu(\lambda )}\, d\lambda \bigg| &\les t^{-\frac{d+1}{2}-\sigma} \\
\label{eq:disp_2wave} \sup_{\xi >0>\xi
'}\bigg|\int_{-\infty}^{\infty}e^{\pm it\lambda}\lambda\chi (\lambda
)\chi_{[|\xi\lambda |>1,|\xi '\lambda
|>1]}(\langle\xi\rangle\langle\xi
'\rangle)^{-\frac{d}{2}-\sigma}\frac{f_{+,\nu}(\xi,\lambda )f_{-,\nu}(\xi ',\lambda
)}{W_\nu(\lambda )}\, d\lambda \bigg| &\les t^{-\frac{d}{2}-\sigma}
\end{align}
\end{lemma}

\begin{proof}
Writing \begin{equation}\label{eq:m+-}
f_{+,\nu}(\xi,\lambda )=e^{i\xi\lambda}m_{+,\nu}(\xi,\lambda
),\qquad f_{-,\nu}(\xi,\lambda )=e^{-i\xi\lambda}m_{-,\nu}(\xi,\lambda
),\end{equation}
 one infers from Lemma~\ref{lem:xge1} that
\[
|m_{+,\nu}(\xi,\lambda)-1|\les \lambda^{-1}\xi^{-1}, \qquad
|\pr_\lambda\, m_{+,\nu}(\xi,\lambda)|\les \lambda^{-2}\xi^{-1}\les
\lambda^{-1}
\]
Indeed, simply set \[ m_{+,\nu}(\xi,\lambda):= (1+\lambda
b_1(\lambda\xi,\lambda))(1+O(\lambda^{-1}\xi^{-1})) \]where $b_1$ is
from that lemma.
Next, we express (\ref{eq:disp_2}) in the form
\begin{equation}\label{eq:47} \bigg|\int_{-\infty}^{\infty}e^{it\phi
(\lambda;\xi,\xi' )}a_\nu(\lambda;\xi,\xi' )\, d\lambda \bigg|\les
t^{-\frac{d+1}{2}-\sigma}
\end{equation} where $\xi
> 0>\xi '$ are fixed, $\phi (\lambda;\xi,\xi'
):=\lambda^2+\frac{\lambda}{t}(\xi -\xi ')$, and
$$a_\nu(\lambda;\xi,\xi' ):=\lambda\chi (\lambda )\chi_{[|\xi\lambda |>1,|\xi '\lambda |>1]}
(\langle\xi\rangle\langle\xi '\rangle
)^{-\frac{d}{2}-\sigma}\frac{m_{+,\nu}(\xi,\lambda )m_{-,\nu}(\xi ',\lambda
)}{W_\nu(\lambda )}$$ Denote the critical point of~$\phi$ by
$\lambda_0:=-\frac{\xi -\xi '}{2t}$.
 By
Proposition~\ref{prop:wronski}, for small $|\lambda |$
$$\Big|\frac{\lambda}{W_\nu(\lambda )}\Big|\les \lambda^{2\nu},\qquad
\Big|\Big(\frac{\lambda}{W_\nu(\lambda )}\Big)'\Big|\les
\lambda^{2\nu-1}$$ Hence,
\begin{align}\label{eq:48} |a_\nu(\lambda;\xi,\xi' )|
&\les \lambda^{2\nu}(\langle\xi\rangle\langle\xi '\rangle
)^{-\frac{d}{2}-\sigma }\chi (\lambda )\chi_{[|\xi\lambda |>1,|\xi '\lambda
|>1]}\\ \label{eq:49} |\pr_\lambda^\ell a_\nu(\lambda;\xi,\xi' )| &\les
\lambda^{2\nu-\ell}(\langle\xi\rangle\langle\xi '\rangle
)^{-\frac{d}{2}-\sigma }\chi (\lambda )\chi_{[|\xi\lambda |>1,|\xi '\lambda
|>1]}
\end{align}
for all $\ell\ge1$. We will need to consider three cases in order to
prove (\ref{eq:disp_2}) via (\ref{46}), depending on where
$\lambda_0$ falls relative to the support of $a$.

\smallskip
\underline{Case 1:} $|\lambda_0 |\les 1$, $|\lambda_0|\gtr |\xi
|^{-1}+|\xi '|^{-1}$.

\smallskip
Note that the second inequality here implies that
$$\frac{\xi +|\xi '|}{t}\gtr\frac{\xi +|\xi '|}{\xi |\xi '|} \text{\ \ or\ \ }1\gtr\frac{t}{\xi |\xi '|}.$$
Furthermore, we remark that $a\equiv 0$ unless $\xi\gtr 1$ and $|\xi
'|\gtr 1$.
Starting with the first integral on the right-hand side
of (\ref{46}) we conclude from (\ref{eq:48}) that
\begin{equation}\label{eq:fall1}
\delta^2\int\frac{|a_\nu(\lambda;\xi,\xi' )|}{|\lambda
-\lambda_0|^2+\delta^2}\;d\lambda \les(\langle\xi\rangle\langle\xi
'\rangle )^{- \frac{d}{2}-\sigma } t^{-1/2}\les  t^{-\frac{d+1}{2}-\sigma}
\end{equation}
For the second integral in~\eqref{46} we obtain from~(\ref{eq:49})
that
\begin{align}
 \delta^2 \int\limits_{|\lambda -\lambda_0|>\delta}\frac{|\pr_\lambda\, a_\nu(\lambda;\xi,\xi'
)|}{|\lambda -\lambda_0|}\, d\lambda &\les
(\langle\xi\rangle\langle\xi '\rangle )^{- \frac{d}{2}-\sigma  } t^{-\frac12} \les t^{-\frac{d+1}{2}-\sigma}
 \label{eq:fall1_der}
\end{align}

\smallskip
\underline{Case 2:} $|\lambda_0 |\les 1$, $|\lambda_0|\ll
\langle\xi\rangle^{-1}+\langle\xi '\rangle^{-1}$.

In this case, we do not use Lemma~\ref{lemma12'}. Instead we note that on the
support of $a_\nu$, we have
$|\partial_\lambda\phi(\lambda;\xi,\xi')|\sim \lambda$, $|\partial_\lambda^2
\phi(\lambda;\xi,\xi')|\les 1$ and the higher derivatives vanish. Let
\[
\lambda_1:= \max(\xi^{-1},|\xi'|^{-1})
\]
Integrating by parts
 thus yields, for sufficiently large $N$,
\begin{align*}
&\Big|\int e^{it\phi(\lambda;\xi,\xi')} \;
a_\nu(\lambda;\xi,\xi')\,d\lambda \Big| \les t^{-N} \int \Big|
(\pr_\lambda (\pr_\lambda \phi)^{-1})^N
a_\nu(\lambda;\xi,\xi')\Big|\,d\lambda \\
&\les t^{-N} \int_{\lambda_1}^1 \lambda^{2\nu-2N}\, d\lambda (\langle\xi\rangle\langle\xi '\rangle
)^{-\frac{d}{2}-\sigma } \\ &\les
t^{-N} (\min(\xi,|\xi'|))^{-2\nu-1+2N} (\langle\xi\rangle\langle\xi '\rangle
)^{-\frac{d}{2}-\sigma } \les t^{-N} (\langle\xi\rangle\langle\xi '\rangle)^{N-\frac{d+1}{2}-\sigma-\nu}\les t^{-\frac{d+1}{2}-\sigma}
\end{align*}
where we used that $\xi|\xi'|\les t$ which follows from
$\lambda_0\ll \lambda_1$.

\smallskip
\underline{Case 3:} $|\lambda_0 |>> 1$,
$|\lambda_0|\gtr\xi^{-1}+|\xi ' |^{-1}$.

\smallskip
In this case, $|\lambda - \lambda_0 |\sim |\lambda_0|>>1$. Thus,
\begin{equation}\label{eq:fall3}
\begin{aligned}
\delta^2 \int\frac{|a(\lambda )|}{|\lambda -\lambda_0|^2+t^{-1}}\, d\lambda
&\les t^{-2} (\langle\xi\rangle\langle\xi '\rangle
)^{\walze} \les  t^{-\frac{d+4}{2}-\sigma}\\
 \delta^2 \int\limits_{|\lambda
-\lambda_0|>\delta}\frac{|a'(\lambda )|}{|\lambda -\lambda_0|}\;
d\lambda &\les  t^{-1} (\langle\xi\rangle\langle\xi '\rangle
)^{\walze}  \les t^{-\frac{d+2}{2}-\sigma}
\end{aligned}
\end{equation} and \eqref{eq:disp_2} is proved.

For \eqref{eq:disp_2wave} note that by $N$-fold integration by parts,
\[
 \Big| \int_{-\infty}^\infty e^{i\lambda(\pm t+\xi-\xi')}\, a_\nu(\lambda;\xi,\xi')\,d\lambda\Big|\les (1+|t+\xi-\xi'|)^{-2\nu-1} (\langle\xi\rangle\langle\xi '\rangle
)^{-\frac{d}{2}-\sigma}
\]
However, since $2\nu+1\ge \nu+\frac12\ge\frac{d}{2}+\sigma$,  this expression is $\les t^{\walze}$ uniformly in $\xi,\xi'$ as claimed.
\end{proof}

Now we turn to the estimate of the oscillatory integral for the case
$|\xi\lambda|
>1$ and $|\xi '\lambda |<1$.

\begin{lemma}\label{lemma14} Let $0\le\sigma\le \nu-\frac{d-1}{2}$.
For all $t>1$ \begin{align}\label{eq:disp_3} \sup_{\xi >0>\xi
'}\bigg|(\langle\xi\rangle\langle\xi '\rangle
)^{\walze}\int_{-\infty}^{\infty}e^{it\lambda^2}\frac{\lambda\chi
(\lambda )}{W_\nu(\lambda )}\chi_{[|\xi\lambda| >1,|\xi '\lambda
|<1]}f_{+,\nu}(\xi,\lambda )f_{-,\nu}(\xi ',\lambda )\, d\lambda \bigg|
&\les
t^{\walz} \\
\label{eq:disp_3wave} \sup_{\xi >0>\xi
'}\bigg|(\langle\xi\rangle\langle\xi '\rangle
)^{\walze}\int_{-\infty}^{\infty}e^{\pm it\lambda}\frac{\lambda\chi
(\lambda )}{W_\nu(\lambda )}\chi_{[|\xi\lambda| >1,|\xi '\lambda
|<1]}f_{+,\nu}(\xi,\lambda )f_{-,\nu}(\xi ',\lambda )\, d\lambda \bigg|
&\les
t^{\walze}
\end{align} and similarly with $\chi_{[|\xi\lambda |<1,|\xi
'\lambda| > 1]}$.
\end{lemma}

\begin{proof}
As before, we write $f_{+,\nu}(\xi,\lambda
)=e^{i\xi\lambda}m_{+,\nu}(\xi,\lambda )$.  But because of $|\xi
'\lambda |<1$ we use the representation
$$f_{-,\nu}(\xi ',\lambda )=a_{-,\nu}(\lambda )u_{0,\nu}^-(\xi',\lambda )+b_{-,\nu}(\lambda )u_{1,\nu}^-(\xi',\lambda )$$
In particular,
Proposition~\ref{prop:wronski}, Lemma~\ref{lem:zero_en}, and
Corollary~\ref{cor4} yield
\begin{equation} |f_{-,\nu}(\xi ',\lambda )| \les
|\lambda|^{\frac12-\nu}\langle\xi '\rangle^{\frac12+\nu},\qquad
|\partial_{\lambda}f_{-,\nu}(\xi ',\lambda )| \les
|\lambda|^{-\frac12-\nu}\langle\xi '\rangle^{\frac12+\nu}
\label{eq:f-m+}
\end{equation}
provided $|\xi '\lambda |<1$. To obtain~\eqref{eq:disp_3} we apply
(\ref{46}) with \[\phi (\lambda
)=\phi(\lambda;\xi,\xi')=\lambda^2+\frac{\xi}{t}\lambda\] and
$$a(\lambda )=a_\nu(\lambda;\xi,\xi')=\frac{\lambda\chi (\lambda )}{W_\nu(\lambda )}(\langle\xi\rangle\langle\xi '\rangle )^{\walze}
\chi_{[|\xi\lambda| >1,|\xi '\lambda |<1]}m_{+,\nu}(\xi,\lambda
)f_{-,\nu}(\xi ',\lambda ).$$
By Proposition~\ref{prop:wronski} and~\eqref{eq:f-m+},
\begin{align}\label{51} |a(\lambda )| &\les |\lambda
|^{\frac12+\nu}  \la\xi\ra^{\walze} \la\xi'\ra^{\nu-\sigma-\frac{d-1}{2}}   \chi (\lambda )\chi_{[|\xi\lambda|
>1,|\xi '\lambda |<1]} \\ \label{52} |\partial_\lambda^\ell a(\lambda )| &\les
|\lambda
|^{\frac12+\nu-\ell}  \la\xi\ra^{\walze} \la\xi'\ra^{\nu-\sigma-\frac{d-1}{2}}   \chi (\lambda )\chi_{[|\xi\lambda|
>1,|\xi '\lambda |<1]} \qquad\forall\;\ell\ge1
\end{align}
The critical point of the phase is $\lambda_0=-\frac{\xi}{2t}$. As usual, we begin with the true stationary phase
case, i.e., $\lambda_0\in\supp(a)$.

\smallskip
\underline{Case 1:} $|\lambda_0|\les 1$, $|\xi\lambda_0|\gtr 1$, $|\xi'\lambda_0|\les 1$

In this case, $|\xi\lambda_0|\gtr 1$ implies that   $\xi\gtr t^{\frac12}$, whereas
 $|\lambda_0|\les 1$ and $|\xi'\lambda_0|\les 1$ together imply that $\la \xi\ra \la \xi'\ra\les t$.
As a consequence, we remark that $|\lambda_0| \gtr \delta=t^{-\frac12}$. Thus, letting $\chi_\delta$ denote a
smooth cutoff to a neighborhood of size $c\delta$ where $c$ is some small positive constant, we conclude that
\begin{align*}
&\Big| \int e^{it\phi(\lambda)} a(\lambda)\,d\lambda \Big| \les \int\limits_{|\lambda-\lambda_0|<c\delta} |\lambda|^{\frac12+\nu} \,d\lambda
\la\xi\ra^{\walze} \la\xi'\ra^{\nu-\sigma-\frac{d-1}{2}}  + t^{-N} \int \Big|\Big(\partial_\lambda \frac{1}{\phi'}\Big)^N (1-\chi_\delta)a\Big|\,d\lambda \\
&\les  \la\xi\ra^{\walze} \la\xi'\ra^{\nu-\sigma-\frac{d-1}{2}} \Big[ \lambda_0^{\frac12+\nu} \delta + t^{-N} \int\limits_{|\lambda-\lambda_0|>c\delta}
\big(|\lambda-\lambda_0|^{-N}\lambda^{\frac12+\nu-N} + |\lambda-\lambda_0|^{-2N} \lambda^{\frac12+\nu}\big)\,d\lambda\Big]
\end{align*}
Carrying out the integrations one checks that the entire right-hand side is $\les t^{-\frac{d+1}{2}-\sigma}$.

\smallskip
\underline{Case 2:} $|\lambda_0|\les 1$, $|\xi'\lambda_0 |> 1$

Then $\la\xi\ra\la\xi'\ra>t$. Using Lemma~\ref{lemma12'} therefore yields
\begin{align*}
\delta^2 \int\frac{|a(\lambda)|}{|\lambda-\lambda_0|^2+\delta^2}\,d\lambda
 &\les \delta^2 \int_{\xi^{-1}}^{\la\xi'\ra^{-1}} \frac{\lambda^{\frac12+\nu} \la\xi\ra^{\walze} \la\xi'\ra^{\nu-\sigma-\frac{d-1}{2}}}{|\lambda - \lambda_0|^2+\delta^2}\,d\lambda \\
&\les \delta \la\xi'\ra^{-\frac12-\nu} \la\xi\ra^{\walze} \la\xi'\ra^{\nu-\sigma-\frac{d-1}{2}} \\
&\les \delta (\la\xi\ra\la\xi'\ra)^{\walze}
\les t^{\walz}
\end{align*}
Similarly,
\begin{align*}
\delta^2 \int\limits_{|\lambda-\lambda_0|>\delta}\frac{|\partial_\lambda a(\lambda)|}{|\lambda-\lambda_0|}\,d\lambda
 &\les \delta \int_{\xi^{-1}}^{\la\xi'\ra^{-1}} \lambda^{-\frac12+\nu} \la\xi\ra^{\walze} \la\xi'\ra^{\nu-\sigma-\frac{d-1}{2}} \,d\lambda \\
&\les \delta \la\xi'\ra^{-\frac12-\nu} \la\xi\ra^{\walze} \la\xi'\ra^{\nu-\sigma-\frac{d-1}{2}}
\les t^{\walz}
\end{align*}
as before.

\smallskip
\underline{Case 3:} $|\lambda_0|\les 1$, $|\xi\lambda_0 |\ll 1$

\smallskip In this case we integrate by parts without using Lemma~\ref{lemma12'}.
 As in the previous lemma, we use that
$|\phi'(\lambda)|\sim\lambda$ on the support of~$a$.  In view of \eqref{51} and~\eqref{52},
\begin{equation}\label{eq:intbyparts}\begin{aligned}
  \Big| \int_{-\infty}^\infty e^{it\phi(\lambda)} a(\lambda)\,
  d\lambda \Big| &\les  t^{-N} \int_{-\infty}^\infty  |
  (\partial_\lambda \phi'(\lambda)^{-1})^N a(\lambda)|\,d\lambda \\
  &\les t^{-N} \la \xi\ra^{\walze} \la\xi'\ra^{\nu-\sigma-\frac{d-1}{2}} \int_{\xi^{-1}}^\infty
  \lambda^{\frac12+\nu-2N}\,d\lambda \\
  & \les t^{-N}
  \xi^{2N-1-d-2\sigma}\les t^{\walz}
\end{aligned}
\end{equation}
where we used that $|\xi'|\le\xi\les t^{\frac12}$.

\smallskip
\underline{Case 4:} $|\lambda_0|>>1$.

\smallskip Here we use Lemma~\ref{lemma12'}.
In view of \eqref{51} and~\eqref{52},
\begin{align*}
\delta^2 \int\frac{|a(\lambda )|}{|\lambda -\lambda_0|^2+t^{-1}}\; d\lambda &\les t^{-2} \la\xi\ra^{\walze} \les t^{-\frac{d+4}{2}-\sigma} \\
\delta^2  \int_{|\lambda -\lambda_0|>\delta}\frac{|a'(\lambda )|}{|\lambda -\lambda_0|}\,d\lambda &\les t^{-1} \la\xi\ra^{\walze} \les t^{-\frac{d+2}{2}-\sigma}
\end{align*}
This proves (\ref{eq:disp_3}).

\medskip
\noindent
For \eqref{eq:disp_3wave}, note that
\[
 \Big| \int_{-\infty}^\infty e^{i\lambda(\pm t+\xi-\xi')}\, a_\nu(\lambda;\xi,\xi')\,d\lambda\Big|\les (1+|\pm t+\xi-\xi'|)^{-\nu-\frac32} \langle\xi\rangle^{-\frac{d}{2}-\sigma}
\]
However, since $\nu+\frac32\ge \nu+\frac12\ge\frac{d}{2}+\sigma$,  this expression is $\les t^{\walze}$ uniformly in $\xi,\xi'$ as claimed.

\noindent
The other case $\chi_{[|\xi\lambda
|<1,\xi '\lambda <-1]}$ is treated in an analogous fashion.
\end{proof}

The remaining cases for the small energy contributions are $\xi
>\xi '>|\lambda |^{-1}$ and $\xi '<\xi <-|\lambda |^{-1}$.  By
symmetry it will suffice to treat the former case.  As usual, we
need to consider reflection and transmission coefficients, therefore
we write \begin{equation}\label{54} f_{-,\nu}(\xi,\lambda
)=\alpha_{-,\nu}(\lambda )f_{+,\nu}(\xi,\lambda )+\beta_{-,\nu}(\lambda
)\overline{f_{+,\nu}(\xi,\lambda )}. \end{equation} Then, with
$W_\nu(\lambda )=W(f_{-,\nu}(\cdot,\lambda ),f_{+,\nu}(\cdot,\lambda ))$,
$$W_\nu(\lambda )=-\beta_-(\lambda )W(f_{+,\nu}(\cdot,\lambda ),\overline{f_{+,\nu}(\cdot,\lambda )})=2i\lambda\beta_{-,\nu}(\lambda )$$
and
\begin{align*}
\wt W_\nu(\lambda):= W(f_{-,\nu}(\cdot,\lambda
),\overline{f_{+,\nu}(\cdot,\lambda )})&=\alpha_{-,\nu}(\lambda
)W(f_{+,\nu}(\cdot,\lambda ), \overline{f_{+,\nu}(\cdot,\lambda )})
=-2i\lambda\alpha_{-,\nu}(\lambda )
\end{align*}
Therefore,
\[
\lambda\frac{\beta_{-,\nu}(\lambda )}{W_\nu(\lambda)} =
-\frac{1}{2i},\qquad \lambda\frac{\alpha_{-,\nu}(\lambda
)}{W_\nu(\lambda)} = \frac{\wt
W_\nu(\lambda)}{W_\nu(\lambda)}=\const+O(\lambda^\eps)
\]
as can be seen from Proposition~\ref{prop:wronski}. The
$O(\lambda^\eps)$ term is complex-valued and behaves like a symbol.

\begin{lemma}\label{lemma15}
For any $t>1$ and $0\le\sigma\le\nu-\frac{d-1}{2}$,
\begin{align}\label{eq:disp_4} \sup_{\xi >\xi
'>0}\bigg|(\langle\xi\rangle\langle\xi '\rangle )^{\walze}\int
e^{it\lambda^2}\frac{\lambda\chi (\lambda )}{W_\nu(\lambda )}\chi_{[|\xi
'\lambda| >1]}f_{+,\nu}(\xi,\lambda )f_{-,\nu}(\xi ',\lambda )\, d\lambda
\bigg| &\les t^{\walz} \\
\label{eq:disp_4wave} \sup_{\xi >\xi
'>0}\bigg|(\langle\xi\rangle\langle\xi '\rangle )^{\walze}\int
e^{\pm it\lambda}\frac{\lambda\chi (\lambda )}{W_\nu(\lambda )}\chi_{[|\xi
'\lambda| >1]}f_{+,\nu}(\xi,\lambda )f_{-,\nu}(\xi ',\lambda )\, d\lambda
\bigg| &\les t^{\walze}
\end{align} and similarly for $\sup_{\xi ' <\xi
<0}$ and $\chi_{[|\xi\lambda |>1]}$.
\end{lemma}

\begin{proof}
Using (\ref{54}), we reduce (\ref{eq:disp_4}) to two estimates, see~\eqref{eq:m+-}:
\begin{align}\label{58}
\sup_{\xi >\xi '>0} (\langle\xi\rangle\langle\xi
'\rangle)^{\walze}\bigg|\int e^{it\lambda^2} e^{i\lambda (\xi +\xi
')} {\chi_{[\xi '|\lambda |>1]}}\chi(\lambda)\, m_+(\xi,\lambda
)m_+(\xi
',\lambda )\, d\lambda \bigg| &\les t^{\walz}  \\
\label{59} \sup_{\xi >\xi '>0}(\langle\xi\rangle\langle\xi '\rangle
)^{\walze} \bigg| \int e^{it\lambda^2}e^{i\lambda (\xi -\xi
')} O(1) \chi_{[\xi '|\lambda |>1]}\chi(\lambda) m_+(\xi,\lambda
)\overline{m_+(\xi ',\lambda )}\, d\lambda \bigg| &\les t^{\walz}
\end{align}
We apply (\ref{46}) to (\ref{58}) with fixed $\xi >\xi '>0$ and
\begin{align*}
\phi (\lambda )& :=\lambda^2+\frac{\lambda}{t}(\xi +\xi '),\\
a(\lambda )& :=(\langle\xi\rangle\langle\xi '\rangle
)^{\walze}\chi_{[\xi '|\lambda |>1]}\chi(\lambda
)m_+(\xi,\lambda )m_+(\xi ',\lambda ).
\end{align*}
Then \begin{align}\label{60} |a(\lambda )|
&\les(\langle\xi\rangle\langle\xi '\rangle )^{\walze}\chi
(\lambda )\chi_{[\xi '|\lambda |>1]} \\
\label{61} |\partial_\lambda^\ell a(\lambda )| &\les |\lambda
|^{-\ell}(\langle\xi\rangle\langle\xi '\rangle )^{\walze}\chi
(\lambda )\chi_{[\xi '|\lambda |>1]} \qquad\forall\;\ell\ge1
\end{align}

\smallskip
\underline{Case 1:} Suppose $|\lambda_0|\les 1$ and $|\xi '\lambda_0
|>1$, where $\lambda_0=-\frac{\xi +\xi '}{2t}$.  Note $\xi >\xi
'\gtr 1$.

\smallskip Then
\begin{align*}
\delta^2\int\frac{|a(\lambda )|}{|\lambda -\lambda_0|^2+t^{-1}}\, d\lambda
&\les \delta^2 (\langle\xi\rangle\langle\xi '\rangle
)^{\walze}\int\frac{d\lambda}{|\lambda -\lambda_0|^2+t^{-1}} \\
&\les t^{-\frac12} t^{\walze}
\end{align*}
since $|\xi '\lambda_0|\sim\frac{\xi\xi '}{t}>1$. As for the
derivative term in (\ref{46}), we infer from (\ref{61}) that
\begin{equation}\label{62} \delta^2\int_{|\lambda -\lambda_0|>\delta}\frac{|a'(\lambda
)|}{|\lambda -\lambda_0 |}\,d\lambda\les \delta^2 (\langle\xi\rangle\langle\xi
'\rangle )^{\walze}\int_{|\lambda
-\lambda_0|>\delta}\frac{d\lambda}{|\lambda ||\lambda -\lambda_0
|}\chi_{[|\lambda\xi '|>1]} \end{equation}
We need to
distinguish between $|\lambda -\lambda_0|>\frac{1}{10}|\lambda_0|$
and $|\lambda -\lambda_0|<\frac{1}{10}|\lambda_0|$.  Thus,
\begin{align*}
(\ref{62})&\les \delta^2 (\langle\xi\rangle\langle\xi '\rangle
)^{\walze}\Big[ \int_{1/\xi
'}^{\infty}\frac{d\lambda}{\lambda^2}+
|\lambda_0|^{-1}\log \big(t^{1/2}|\lambda_0|\big)\Big]\\
& \les \delta^2 (\langle\xi\rangle\langle\xi '\rangle
)^{\walze} \Big[ \xi' + t\xi^{-1}\log(\xi t^{-\frac12})\Big] \les \delta^2  (\langle\xi\rangle\langle\xi '\rangle
)^{\walze+\frac12}\\
&\les t^{\walz}
\end{align*}
where we used $\xi>\xi'\gtr1$, $\xi\xi'\gtr t$, and $\xi^2>t$.

\smallskip
\underline{Case 2:} $|\lambda_0|\les 1$, $|\lambda_0|\ll
\frac{1}{\xi '}$.

\smallskip Then $|\phi'(\lambda)|\sim|\lambda |$ on the support
of $a(\lambda)$.  Hence, integration by parts yields
\begin{align*}
&\Big|\int e^{it\phi(\lambda)} a(\lambda)\;d\lambda\Big| \les t^{-N} \int \Big|(\partial_\lambda \phi'(\lambda)^{-1})^N a(\lambda)\Big|\,d\lambda \\
&\les  t^{-N} (\langle\xi\rangle \langle\xi'\rangle
)^{\walze} \int\limits_{\frac{1}{\xi'}}^\infty \lambda^{-2N}\,d\lambda  \les t^{-N} (\xi')^{2N-d-1-2\sigma}\les t^{\walz}
\end{align*}
where we used that $(\xi')^2<\xi\xi'\les t$.

\smallskip
\underline{Case 3:} $|\lambda_0|>>1$, $|\lambda_0|\gtr\frac{1}{\xi
'}$.

\smallskip
Then $|\lambda -\lambda_0|\sim |\lambda_0|$ on $\supp (a)$ and $\xi>t$.
Using Lemma~\ref{lemma12'} yields
\begin{align*}
\delta^2\int\frac{|a(\lambda )|}{|\lambda -\lambda_0|^2+t^{-1}}\, d\lambda\les
t^{-2} (\langle\xi\rangle \langle\xi'\rangle
)^{\walze} \les t^{-\frac{d+4}{2}-\sigma}
\end{align*}
as well as
\begin{align*}
\delta^2 \int_{|\lambda -\lambda_0|>\delta}\frac{|a'(\lambda )|} {|\lambda
-\lambda_0|}\, d\lambda&\les \delta^2 (\langle\xi\rangle\langle\xi '\rangle
)^{\walze} |\lambda_0|^{-1}\int_{\frac{1}{\langle\xi
'\rangle}}^{1}\frac{d\lambda}{|\lambda |}\\
&\les \delta^2 (\langle\xi\rangle\langle\xi '\rangle
)^{\walze} \frac{t}{\xi}\log\la\xi'\ra \les \delta^2 (\langle\xi\rangle\langle\xi '\rangle
)^{\walze} t\la\xi\ra^{-\frac12}\\
&\les \delta (\langle\xi\rangle\langle\xi '\rangle
)^{\walze} \les t^{\walz}.
\end{align*}
This concludes the proof of (\ref{58}).  For~\eqref{59} we argue analogously, with $\lambda_0=-\frac{\xi-\xi'}{2t}$.
Case~1 above applies without major changes since we again have $\xi\xi'\gtr t$. In Case~2, however, we cannot guarantee
that $\xi\xi\les t$ as before (since we could have $\xi=\xi'$ and $\lambda_0=0$, say). However, if $\xi\xi'\gtr t$, then
the calculations of Case~1 yield the desired conclusion. Finally, Case~3 is the same and \eqref{59} is proved.

\medskip\noindent For \eqref{eq:disp_4wave}  we first note the bound, obtained by repeated integration by parts using~\eqref{61}
\[
 \Big| \int_{-\infty}^\infty e^{i\lambda(\pm t+\xi-\xi')}\, a_\nu(\lambda;\xi,\xi')\,d\lambda\Big|\les (1+|\pm t+\xi-\xi'|)^{-N} (\langle\xi\rangle \langle\xi'\rangle)^{-\frac{d}{2}-\sigma} |\xi'|^{N-1}
\]
valid for all integers $N\ge1$.  Now let us choose $N$ such that $N\ge\frac{d}{2}+\sigma\ge N-1$. Then  this bound is dominated by $t^{\walze}$ uniformly in~$\xi,\xi'$ as claimed.

The case of $\xi '<\xi <0$, $|\xi\lambda |>1$
is treated analogously.
\end{proof}

We are done with the contributions of small $\lambda$ to our main
oscillatory integral. To conclude the proof of Theorem~\ref{thm1} it
suffices to prove the following statement about the contributions
from ``large'' energies.

\begin{lemma}\label{lemma16}
For all $t>0$, \begin{align}\label{eq:disp_large} \sup_{\xi >\xi
'}\bigg|(\langle\xi\rangle\langle\xi '\rangle
)^{\walze}\int_{-\infty}^{\infty}e^{it\lambda^2}\frac{\lambda
(1-\chi )(\lambda )}{W_\nu(\lambda )}f_{+,\nu}(\xi,\lambda )f_{-,\nu}(\xi
',\lambda )\, d\lambda \bigg| &\les  t^{\walz}
 \end{align}
 where $\sigma\ge0$ is arbitrary.
\end{lemma}

\begin{proof}
We observed above, see (\ref{54}), that $W_\nu(\lambda
)=-2i\lambda\beta_{-,\nu}(\lambda )$. Since $|\beta_{-,\nu}(\lambda
)|\geq 1$, this implies that $|W_\nu(\lambda )|\geq 2|\lambda |$. In
particular, $W_\nu(\lambda )\neq 0$ for every $\lambda\neq 0$. We
shall write
\begin{equation}\nonumber
f_{+,\nu}(\xi,\lambda) = e^{i\xi\lambda} m_{+,\nu}(\xi,\lambda),\qquad
f_{-,\nu}(\xi,\lambda) = e^{-i\xi\lambda} m_{-,\nu}(\xi,\lambda)
\end{equation}
The functions $m_{\pm,\nu}(\xi,\lambda)$ satisfy the Volterra equation
\begin{equation}\label{65}
m_{+,\nu}(\xi,\lambda
)=1+\int_{\xi}^{\infty}\frac{1-e^{-2i(\tilde{\xi} -\xi
)\lambda}}{2i\lambda}\wt V_\nu(\xi)(\tilde{\xi}
)m_{+,\nu}(\tilde{\xi},\lambda )d\tilde{\xi}
\end{equation}
where the potential
$\wt V_\nu$
satisfies
$$\bigg|\frac{d^\ell}{d\xi^\ell}\wt V_\nu(\xi )\bigg|\les\langle\xi\rangle^{-2-\ell},\quad\forall\; \ell\geq 0.$$
From (\ref{65}), for any $\xi\geq 0$
$$m_{+,\nu}(\xi,\lambda )=1+O(\lambda^{-1}\langle\xi\rangle^{-1})$$
Moreover, see~\cite{SSS},
\[
|\partial_{\lambda}^k \partial_{\xi}^\ell m_{+,\nu}( \xi,\lambda )| \le
C_{k,\ell}\, \lambda^{-1-k}\langle\xi\rangle^{-1-\ell}
\quad\forall\; k+\ell>0
\]
In~\cite{SSS} this is proved for $k+\ell\le2$, but the proof there
extends inductively to higher orders. As a corollary, we obtain
(take $\xi =0$)
\begin{align*}
W_\nu(\lambda )&=W(f_{-,\nu}(\cdot,\lambda ),f_{+,\nu}(\cdot,\lambda ))\\
&=m_{+,\nu}(\xi,\lambda )[m'_{-,\nu}(\xi,\lambda )-i\lambda m_{-,\nu}(\xi,\lambda
)]-m_{-,\nu}(\xi,\lambda )
[m'_{+,\nu}(\xi,\lambda )+i\lambda m_{+,\nu}(\xi,\lambda )]\\
&=-2i\lambda (1+O(\lambda^{-1}))+O(\lambda^{-1}) =-2i\lambda+O(1)
\end{align*}
with derivatives $(\lambda/W_\nu(\lambda
))^{(\ell)}=O(\lambda^{-1-\ell})$ as $|\lambda |\to\infty$.

In order to prove (\ref{eq:disp_large}), we will need to distinguish
the cases $\xi
>0>\xi '$, $\xi
>\xi '>0$, as well as $0>\xi >\xi '$.  By symmetry, it will suffice to consider
the first two.

\smallskip
\underline{Case 1:} $\xi >0>\xi '$.

\smallskip\noindent
In this case we need to prove that
\begin{align}\label{64} \sup_{\xi
>0>\xi '}\bigg|(\langle\xi\rangle\langle\xi '\rangle
)^{-\frac{d}{2}-\sigma}\int e^{it[\lambda^2+\frac{\xi -\xi
'}{t}\lambda]}\frac{\lambda (1-\chi )(\lambda )}{W_\nu(\lambda )}\,
m_{+,\nu}(\xi,\lambda )m_{-,\nu}(\xi ',\lambda )\, d\lambda \bigg| &\les
t^{\walz}
\end{align}
Let $\phi (\lambda ):=\lambda^2+\frac{\xi -\xi '}{t}\lambda$ and
$$a(\lambda ):= (\langle\xi\rangle\langle\xi '\rangle )^{-\frac{d}{2}-\sigma}
\frac{\lambda (1-\chi )(\lambda )}{W_\nu(\lambda )}m_{+,\nu}(\xi,\lambda
)m_{-,\nu}(\xi ',\lambda ).$$ The phase $\phi$ has critical
point~$\lambda_0=-\frac{\xi -\xi '}{2t}$. If $|\lambda_0|\ll 1$,
then $|\phi'(\lambda)|\sim |\lambda|$ on the support of~$a$.
Therefore, integrating by parts repeatedly yields that the
oscillatory integral in~\eqref{64} decays like $t^{-N}$ for
all~$N\ge1$. It is important to note that already a single
derivative in~$\lambda$ renders $a(\lambda)$ of
size~$O(\lambda^{-2})$ and thus integrable.

For the remainder of Case~1 we can therefore assume that
$|\lambda_0|\gtr 1$ which implies that $\max(\xi,|\xi'|)\gtr t$. In
particular, the factor $(\langle\xi\rangle\langle\xi '\rangle
)^{\walze}$ contributes the decay~$t^{\walze}$.  Applying
Lemma~\ref{lemma12'} yields
\begin{align*}
(\ref{64})&\les \delta^2 \int\frac{|a(\lambda )|}{|\lambda
-\lambda_0|^2+t^{-1}}\,d\lambda
+ \delta^2\int\limits_{|\lambda -\lambda_0|>\delta}\frac{|a'(\lambda )|}{|\lambda -\lambda_0|}\,d\lambda\\
&=: A+B
\end{align*}
By the preceding,
$$A\les \delta \|a\|_{\infty}\les t^{-\frac{1}{2}}(\langle\xi\rangle\langle\xi '\rangle )^{\walze}
\les t^{\walz} $$
\noindent  Next, we estimate $B$. First, from our bounds
on $W_\nu(\lambda )$, $m_{+,\nu}(\xi,\lambda )$ and $m_{-,\nu}(\xi
',\lambda )$ we conclude that
$$|a'(\lambda )|\les (\langle\xi\rangle\langle\xi '\rangle )^{\walze}\chi_{[|\lambda |\gtr 1]}|\lambda |^{-2}.$$
Second, because of  $|\lambda_0|\gtr 1$ we obtain
\begin{align*}
B&\les \delta^2 (\langle\xi\rangle\langle\xi '\rangle )^{\walze}
\int_{\binom{|\lambda -\lambda_0|>\delta}{|\lambda |\gtr 1}}\frac{d\lambda}{|\lambda |^2\, |\lambda -\lambda_0|}\\
&\les \delta (\langle\xi\rangle\langle\xi '\rangle
)^{\walze}\int_{|\lambda |\gtr 1}\frac{d\lambda}
{\lambda^2}\les t^{\walz}.
\end{align*}
This finishes the case $\xi >0>\xi '$.

\smallskip
\underline{Case 2:} $\xi >\xi '>0$

From (\ref{54})
$$f_{-,\nu}(\xi ',\lambda )=\alpha_{-,\nu}(\lambda )f_{+,\nu}(\xi ',\lambda )+\beta_{-,\nu}(\lambda )\overline{f_{+,\nu}(\xi ',\lambda )}$$
where
\begin{align*}
\alpha_{-,\nu}(\lambda )&=\frac{W(f_{-,\nu}(\cdot,\lambda ),\overline{f_{+,\nu}(\cdot,\lambda )})}{-2i\lambda}\\
\beta_{-,\nu}(\lambda )&=\frac{W(f_{+,\nu}(\cdot,\lambda ),f_{-,\nu}(\cdot,\lambda
))}{-2i\lambda}=\frac{W_\nu(\lambda )}{2i\lambda}
\end{align*}
From our large $\lambda$ asymptotics of $W_\nu(\lambda )$ we deduce that
\begin{equation}\label{eq:beta-}\beta_{-,\nu}(\lambda )=1+O(\lambda ^{-1}), \quad \beta_{-,\nu}'(\lambda
)=O(\lambda^{-2}).
\end{equation}
For $\alpha_{-,\nu}(\lambda )$ we calculate, again at
$\xi =0$,
\begin{align*}
W(f_{-,\nu}(\cdot,\lambda ),\overline{f_{+,\nu}(\cdot,\lambda
)})=\,&m_{-,\nu}(\xi,\lambda )(\overline{m}_{+,\nu}'(\xi,\lambda )
-2i\lambda\overline{m}_{+,\nu}(\xi,\lambda ))\\
&-\overline{m}_{+,\nu}(\xi,\lambda )(m'_{-,\nu}(\xi,\lambda )-2i\lambda m_{-,\nu}(\xi,\lambda ))\\
=\,&m_{-,\nu}(\xi,\lambda )\overline{m}'_{+,\nu}(\xi,\lambda)-m_{-,\nu}'(\xi,\lambda )\overline{m}_{+,\nu}(\xi,\lambda )\\
=\,&O(\lambda^{-1})
\end{align*}
so that
\begin{equation}\label{eq:alpha-}\alpha_{-,\nu}(\lambda )=O(\lambda^{-2}),\quad \alpha'_{-,\nu}(\lambda )=O(\lambda^{-3}).
\end{equation}
Thus, we are left with proving the two bounds
\begin{align}\label{70}
\sup_{\xi >\xi
'>0}&\bigg|\int_{-\infty}^{\infty}e^{it\lambda^2}e^{i\lambda (\xi
+\xi ')}\frac{\lambda (1-\chi (\lambda ))}{W_\nu(\lambda
)}\alpha_{-,\nu}(\lambda )\frac{m_{+,\nu}(\xi,\lambda )m_{+,\nu}(\xi ',\lambda
)}{(\langle\xi\rangle\langle\xi '\rangle)^{\frac{d}{2}+\sigma}}\,
d\lambda \bigg|
\les t^{\walz} \\
\label{71} \sup_{\xi >\xi
'>0}&\bigg|\int_{-\infty}^{\infty}e^{it\lambda^2}e^{i\lambda (\xi
-\xi ')}\frac{\lambda (1-\chi (\lambda ))}{W_\nu(\lambda
)}\beta_{-,\nu}(\lambda )\frac{m_{+,\nu}(\xi,\lambda )\overline{m_{+,\nu}(\xi
',\lambda )}}{(\langle\xi\rangle\langle\xi
'\rangle)^{\frac{d}{2}+\sigma}}\, d\lambda \bigg| \les t^{\walz}
\end{align}
for any $t>0$. This, however, follows by means of the exact same
arguments which we use to prove (\ref{64}).  Note that in (\ref{70})
the critical point of the phase is
$$\lambda_0=-\frac{\xi +\xi '}{2t}$$
whereas in (\ref{71}) it is $\lambda_0=-\frac{\xi -\xi '}{2t}.$ In
either case it follows from $|\lambda_0|\gtr 1$ that $\xi\gtr t$.
Hence we can indeed argue as in Case~1. This finishes the proof of
the lemma, and thus also establishes Theorem~\ref{thm1}.
\end{proof}

Now for the wave case. We will tacitly use some elements of the
previous proof.

\begin{lemma}\label{lem:wave_crux}
For all $t>0$, \begin{align}\label{eq:wave_biglam} &
\bigg|\int_{-\infty}^\xi (\langle\xi\rangle\langle\xi '\rangle
)^{\walze}\int_{-\infty}^{\infty}e^{\pm
it\lambda}\frac{\lambda (1-\chi )(\lambda )}{W_\nu(\lambda
)}f_{+,\nu}(\xi,\lambda )f_{-,\nu}(\xi ',\lambda
)\, d\lambda\;  \phi(\xi')\, d\xi'\bigg| \nonumber\\
& \les t^{\walze} \int \big(|\phi(\eta)|+|\phi'(\eta)|\big)\, d\eta.
\end{align}
with a constant that does not depend on $\xi$.
\end{lemma}

\begin{proof}
In order to prove (\ref{eq:wave_biglam}), we will need to
distinguish the cases $\xi >0>\xi '$, $\xi >\xi '>0$, and $0>\xi
>\xi '$.  By symmetry, it will suffice to consider the first two.

\smallskip
\underline{Case 1:} $\xi >0>\xi '$.

\smallskip
Integrating by parts yields
 \begin{align*} &\bigg|(\langle\xi\rangle\langle\xi '\rangle
)^{\walze}\int e^{i\lambda(\pm t+\xi -\xi ')}\frac{\lambda
(1-\chi )(\lambda )}{W_\nu(\lambda )}m_{+,\nu}(\xi,\lambda )m_{-,\nu}(\xi ',\lambda
)\, d\lambda \bigg|\\
&\les (\langle\xi\rangle\langle\xi '\rangle )^{\walze}
|t\pm(\xi-\xi')|^{-N} \les t^{\walze}
\end{align*}
provided $|t\pm(\xi-\xi')|\ge1$. If this fails, then we need to
integrate by parts in~$\xi'$ to remove one factor of~$\lambda$:
since $\lambda e^{-i\xi'\lambda} = i\partial_{\xi'}
e^{-i\xi'\lambda}$, it follows that
\begin{align*}
  & \int_{-\infty}^\xi \langle\xi\rangle^{\walze} \langle\xi '\rangle^{\walze}
  \int e^{i\lambda(\pm t+\xi -\xi ')}\frac{\lambda
(1-\chi )(\lambda )}{W_\nu(\lambda )}m_{+,\nu}(\xi,\lambda )m_{-,\nu}(\xi ',\lambda
)\, d\lambda\, \phi(\xi')\,d\xi' =\\
& i\la \xi\ra^{-d-2\sigma} \int e^{\pm it\lambda} \frac{(1-\chi )(\lambda
)}{W_\nu(\lambda )}m_{+,\nu}(\xi,\lambda )m_{-,\nu}(\xi,\lambda )\, d\lambda\,
\phi(\xi)
\\
 &-i\int_{-\infty}^\xi \langle\xi\rangle^{\walze}
  \int e^{i\lambda(\pm t+\xi -\xi ')}\frac{
(1-\chi )(\lambda )}{W_\nu(\lambda )}m_{+,\nu}(\xi,\lambda )
\partial_{\xi'}\big[\langle\xi '\rangle^{\walze}m_{-,\nu}(\xi ',\lambda )\, \phi(\xi')\big]\,
d\lambda\,d\xi'.
\end{align*}
Denote the two expressions after the equality sign by $A$ and $B$, respectively.
First, exploiting the cancellation due to $W_\nu(-\lambda)=-W_\nu(\lambda)+ O(1)$ as $\lambda\to\infty$, we see that
\[
\sup_{\xi>0>\xi'}\Big|\int e^{it\lambda} \frac{(1-\chi )(\lambda
)}{W_\nu(\lambda )}m_{+,\nu}(\xi,\lambda )m_{-,\nu}(\xi ,\lambda )\, d\lambda\Big|
\les 1 .
\]
Furthermore, integrating by parts in $\lambda$ shows that the left-hand side is in fact $\les t^{-N}$ for any $N$.
Hence,
\[
A\les \la t\ra^{-N} \sup|\phi|\le \la t\ra^{-N} \int (|\phi'(\xi)|+|\phi(\xi)|)\, d\xi
\]
Second, by the same cancellation,
\begin{align*}
B&\les \int (\la\xi\ra\la\xi'\ra)^{\walze} (1+|t\pm (\xi-\xi')|)^{-N}  (|\phi'(\xi')|+|\phi(\xi')|)\, d\xi'\\
&\les t^{\walze}  \int (|\phi'(\xi')|+|\phi(\xi')|)\, d\xi'
\end{align*}
which gives the desired bound as usual.

\smallskip
\underline{Case 2:}  $\xi >\xi '>0$

\medskip

In analogy with \eqref{70} and \eqref{71} we need to consider
\begin{align}\label{70'}
& \int_{-\infty}^{\infty}e^{it\lambda}e^{i\lambda (\xi
+\xi ')}\frac{\lambda (1-\chi (\lambda ))}{W_\nu(\lambda
)}\alpha_{-,\nu}(\lambda )\frac{m_{+,\nu}(\xi,\lambda )m_{+,\nu}(\xi ',\lambda
)}{(\langle\xi\rangle\langle\xi '\rangle)^{\frac{d}{2}+\sigma}}\, d\lambda \\
\label{71'} &\int_{-\infty}^{\infty}e^{it\lambda}e^{i\lambda (\xi
-\xi ')}\frac{\lambda (1-\chi (\lambda ))}{W_\nu(\lambda
)}\beta_{-,\nu}(\lambda )\frac{m_{+,\nu}(\xi,\lambda )\overline{m_{+,\nu}(\xi
',\lambda )}}{(\langle\xi\rangle\langle\xi '\rangle)^{\frac{d}{2}+\sigma}}\,
d\lambda
\end{align}
The integral in \eqref{70'} is $\les \la t\ra^{\walze}$ uniformly in $\xi,\xi'$ due to the decay of~$\alpha_{-,\nu}$, see~\eqref{eq:alpha-}.
On the other hand, the integral in~\eqref{71'} is not a bounded function in $\xi,\xi'$ due to the lack of decay in~$\lambda$, see~\eqref{eq:beta-}.
Thus, we again need to redeem one power of $\lambda$ via a $\xi'$ differentiation, see above.
\end{proof}

\end{document}